\documentclass[a4paper,11pt]{amsart}

\usepackage{amssymb}
\usepackage{amsfonts}
\usepackage{amsmath}
\usepackage{hyperref}
\usepackage[bbgreekl]{mathbbol}

\usepackage{graphicx}
\usepackage[all]{xy}
\usepackage{array}

\usepackage{amsthm}

\newtheorem{Theorem}{Theorem}[section]
\newtheorem{lemme}[Theorem]{Lemma}
\newtheorem{Proposition}[Theorem]{Proposition}
\newtheorem{Corollary}[Theorem]{Corollary}

\theoremstyle{definition}
\newtheorem{Definition}[Theorem]{Definition}
\newtheorem{Notation}[Theorem]{Notation}

\theoremstyle{remark}
\newtheorem*{Remark}{Remark}

\newcommand{\Si}{\mathfrak{S}}
\newcommand{\Ext}{{\mathrm{Ext}}}
\renewcommand{\hom}{{\mathrm{Hom}}}
\newcommand{\Id}{{\mathrm{Id}}}

\title{Universal classes for algebraic groups}
\author{Antoine Touz\'e}
\date{
\today}
\begin{document}

\maketitle
\sloppy
\begin{abstract}
We exhibit cocycles representing certain classes in the cohomology of the 
algebraic group $GL_{n}$ with coefficients in the representation 
$\Gamma^*(\mathfrak{gl}_n^{(1)})$. These classes' existence was anticipated 
by van der Kallen, and they intervene in the proof that reductive 
linear algebraic groups have finitely generated cohomology algebras 
\cite{TvdK}. 
\end{abstract}

Let $\mathbb{k}$ be a field of positive characteristic, let $A$ be a finitely
generated $\mathbb{k}$-algebra, and let $G$ be a reductive linear algebraic
group defined over $\mathbb{k}$ and acting rationally on $A$ by algebra 
automorphisms. Then the rational  cohomology $H^*(G,A)$ is an algebra, 
and one can wonder if it is finitely generated. In degree $0$, the 
finite generation of the subalgebra $A^G=H^0(G,A)$ is part of 
Hilbert's fourteenth problem and was solved positively by the 
work of Nagata \cite{Nagata} and Haboush \cite{Haboush}. 
The finite generation of the whole cohomology algebra remained 
unsolved in general, though much progress had been made in recent 
years \cite{FS,VdK}.

In \cite{VdK}, van der Kallen proved (under some restrictions on the
characteristic of $\mathbb{k}$ which were  removed in \cite{SvdK}) that the
finite generation of $H^*(G,A)$ holds under the following condition: the group $G$ embeds in $GL_{n,\mathbb{k}}$ for some integer $n$, and there exist universal cohomology classes in $ H^*(GL_{n,\mathbb{k}},\Gamma^*(\mathfrak{gl}_n^{(1)}))$ satisfying some divided power algebra relations. He proved  the existence of these universal cohomology classes for $n=2$ \cite[Th 4.4]{VdK}, and for $n=3$ in characteristic $p=2$ \cite{VdK2}.

Later, van der Kallen mentioned that cohomological finite generation holds under a weaker condition, namely the existence of the so-called ``lifted universal cohomology classes''. Our main result is the existence of these lifted classes:
\begin{Theorem}\label{thm-lifted-classes}
Let $\mathbb{k}$ be a field of positive characteristic and let $n>1$ be an integer. There are cohomology classes $c[d]\in H^{2d}(GL_{n,\mathbb{k}},\Gamma^d(\mathfrak{gl}_n^{(1)}))$ such that~:
\begin{enumerate}
\item $c[1]\in H^2(GL_{n,\mathbb{k}},\mathfrak{gl}_n^{(1)})$ is non zero.
\item Let $d\ge 1$ and let $\Delta_{(1,\dots,1)} : \Gamma^d(\mathfrak{gl}_n^{(1)})\to (\mathfrak{gl}_n^{(1)})^{\otimes d}$ be the map induced by the diagonal $\Gamma^d\to \otimes^d$. Then $\Delta_{(1,\dots,1)\,*}c[d]=c[1]^{\cup d}$.
\end{enumerate}
\end{Theorem}

This actually reproves and significantly 
extends the famous theorem of
Friedlander and Suslin \cite[Thm 1.2]{FS}
on the existence of nonzero universal
classes $e_r$. Indeed, as indicated in 
\cite[Lemma 6.2]{TvdK}, we may draw the following corollary from theorem 
\ref{thm-lifted-classes}:
\begin{Corollary} (\cite[Thm 1.2]{FS}) Let $\mathbb{k}$ be a field of positive
characteristic. For any $n>1$ there exist classes
$$e_r\in H^{2p^{r-1}}(GL_{n,\mathbb{k}}, \mathfrak{gl}_n^{(r)})$$
which restrict nontrivially to $H^{2p^{r-1}}((GL_{n,\mathbb{k}})_{(1)},
\mathfrak{gl}_n^{(r)})$.
 \end{Corollary}

The proof of theorem \ref{thm-lifted-classes} may be summarized as follows. In section \ref{sec-cohom-bif}, we remark that theorem \ref{thm-lifted-classes} reduces to a stable cohomology statement, that is, it suffices to prove it for large values of $n$. Bifunctor cohomology \cite{FF} gives access to the stable rational cohomology of $GL_{n,\mathbb{k}}$, and we translate theorem \ref{thm-lifted-classes} in terms of (strict polynomial) bifunctors. More specifically, we show that theorem \ref{thm-lifted-classes} reduces to theorem \ref{thm-lifted-classes-cohom-bif}, that is, to the computation of some classes in the cohomology of the strict polynomial bifunctors $\Gamma^d(gl^{(1)})$.

The proof of theorem \ref{thm-lifted-classes-cohom-bif} is given in section
\ref{sec-proof-existence-lifted-classes}. We use explicit coresolutions of the
bifunctors $\Gamma^d(gl^{(1)})$ to compute cocycles representing the classes $c[d]$.

Section \ref{sec-explicit} is devoted to building the explicit coresolutions of the $\Gamma^d(gl^{(1)})$. We use the following strategy.
First, the $\Gamma^d(gl^{(1)})$ are bifunctors of the form $F(gl)$, obtained by precomposing a  functor $F$ by the bifunctor $gl(-_1,-_2):=\hom_{\mathbb{k}}(-_1,-_2)$.  We remark in section \ref{sec-cohom-Fgl} that the cohomology of this kind of bifunctor may be computed \emph{via} acyclic coresolutions obtained by precomposing an injective coresolution of $F$ by the bifunctor $gl$. Thus, our seeking of the explicit coresolutions of the bifunctors $\Gamma^d(gl^{(1)})$ reduces to the (combinatorially easier) seeking of injective coresolutions of the functors $\Gamma^d(I^{(1)})$ obtained by precomposing the functors $\Gamma^d$ by the Frobenius twist $I^{(1)}$. 
Second, we define in section \ref{sec-twist-comp-res} a class of injective coresolutions of strict polynomial functors called ``twist compatible coresolutions". These coresolutions enjoy the following nice property: we may use  an injective twist compatible coresolution $J_F$ of $F$ to build an explicit injective coresolution of the functor $F(I^{(1)})$.
Third, we build injective twist compatible coresolutions of the functors $\Gamma^d$ in section \ref{sec-twist-comp-res-Gammagl}.

When the characteristic $p$ is odd, the combinatorics of the Frobenius twist bring the notion of $p$-complex into play. 
Section \ref{sec-p-complexes} contains the recollections about $p$-complexes needed in section \ref{sec-twist-comp-res}, as well as a new result (proposition \ref{prop-comparaison-prod-tens}) which is the key point to identify cup products in section \ref{sec-proof-existence-lifted-classes}.
\tableofcontents

\section{Reduction to a bifunctor cohomology computation}\label{sec-cohom-bif}
Let $\mathbb{k}$ be a field of positive characteristic. Using functors is a classical way to build representations of (ordinary) groups: if a group $G$ acts linearly on a $\mathbb{k}$-vector space $V$ and if $F$ is a functor from $\mathbb{k}$-vector spaces to $\mathbb{k}$-vector spaces, then functoriality endows $F(V)$ with an action of $G$. To do the same in the framework of algebraic groups and rational representations, one has to use an algebraic modification of the notion of functor, namely the ``strict polynomial functors'' introduced in \cite{FS}. Classical functors such as the tensor products $\otimes^d$, the divided powers $\Gamma^d$, the symmetric powers $S^d$ or the Frobenius twist $I^{(1)}$ are strict polynomial functors. Now the representations we are interested in are not given by strict polynomial functors but by strict polynomial \emph{bi}functors, contravariant in the first variable and covariant in the second one, as introduced in \cite{FF}. Examples of such bifunctors are the bifunctor $gl(-_1,-_2):=\hom_{\mathbb{k}}(-_1,-_2)$ or postcompositions of $gl$ by strict polynomial functors, such as $gl^{(1)}= I^{(1)}\circ gl$ or $\Gamma^d(gl^{(1)})=\Gamma^d\circ I^{(1)} \circ gl$.

Thus, if $B$ is a strict polynomial bifunctor over $\mathbb{k}$, then for all integer $n$ the $\mathbb{k}$-vector space $B(\mathbb{k}^n,\mathbb{k}^n)$ is endowed with an action of $GL_{n,\mathbb{k}}$. Such representations, with $n$ big enough are called ``stable'' representations of $GL_{n,\mathbb{k}}$. Let's specify what ``big enough'' means.
The polynomial nature of strict polynomial bifunctors endows them with a notion of bidegree. For example, $gl$ is a homogeneous bifunctor of bidegree $(1,1)$ and if $F$ is a  homogeneous strict polynomial functor of degree $d$ then the composite $F(gl)$ is a homogeneous bifunctor of bidegree $(d,d)$. We denote by $\mathcal{P}^d_{e,\mathbb{k}}$ the abelian category (quite strangely denoted by $\mathcal{P}_d^{\mathrm{op}}\times \mathcal{P}_e$ in \cite{FF}) of homogeneous strict  polynomial bifunctors of bidegree $(d,e)$, defined over $\mathbb{k}$, contravariant in the first variable and covariant in the second one.  A stable representation is a representation of the form $B(\mathbb{k}^n,\mathbb{k}^n)$, with $B\in\mathcal{P}^d_{e,\mathbb{k}}$ and $n\ge\max(d,e)$.

$\Ext$-computations in the strict polynomial bifunctor categories $\mathcal{P}^d_{e,\mathbb{k}}$ give access \cite[Th 1.5]{FF} to the rational cohomology of $GL_{n,\mathbb{k}}$ with coefficients in stable representations. In this section we first prove a strengthening of this result, and next we use it to translate theorem \ref{thm-lifted-classes} into the world of strict polynomial bifunctors.

\subsection{Bifunctor cohomology and stable rational cohomology}\label{subsec-cohom-bif-rat-cohom}

Let $\mathbb{k}$ be a field of positive characteristic. The rational cohomology of $GL_{n,\mathbb{k}} $ with coefficients in the rational $GL_{n,\mathbb{k}}$-module $M$ is defined as the extension groups in the category $GL_{n,\mathbb{k}}\text{-mod}$ of rational $GL_{n,\mathbb{k}}$-modules : $$H^*(GL_{n,\mathbb{k}},M):=\Ext^*_{GL_{n,\mathbb{k}}\text{-mod}}(\mathbb{k},M)\,,$$
where $\mathbb{k}$ has a trivial $GL_{n,\mathbb{k}}$-module structure. Similarly, the cohomology of a bifunctor $B$ homogeneous of bidegree $(d,d)$ is defined as the extension groups:
$$H^*_{\mathcal{P}}(B):=\Ext^*_{\mathcal{P}_{d,\mathbb{k}}^d}(\Gamma^d (gl), B)\;.$$
We already know \cite[Th 1.5]{FF} that $H^*_{\mathcal{P}}(B)$ is isomorphic to the stable rational cohomology of $GL_{n,\mathbb{k}}$ with coefficients in $B(\mathbb{k}^n,\mathbb{k}^n)$. We give a more explicit description of this isomorphism which shows its compatibility with cup products.

If $B$ is a bifunctor and if $E$ is a $m$-fold extension of $\Gamma^d (gl)$ by $B$, we may evaluate it on the pair $(\mathbb{k}^n,\mathbb{k}^n)$ and pull it back by the $GL_{n,\mathbb{k}}$-equivariant morphism $\mathbb{k}\to \Gamma^d(gl(\mathbb{k}^n,\mathbb{k}^n))$, 
$x\mapsto x (\Id_{\mathbb{k}^n})^{\otimes d}$.
This defines a map (natural in $B$)
$$\phi_{B,n}:H^*_{\mathcal{P}}(B)\to H^*(GL_{n,\mathbb{k}},B(\mathbb{k}^n,\mathbb{k}^n))\;.$$
This map is an isomorphism in the stable range. More specifically:
\begin{lemme}\label{lm-cohom-bif-computes-cohom-rat}
Let $\mathbb{k}$ be a field of characteristic $p>0$ and $B\in \mathcal{P}_{d,\mathbb{k}}^d$ a homogeneous functor of bidegree $(d,d)$. If $n\ge d$ then $\phi_{B,n}$ is an isomorphism.  
\end{lemme}

\begin{proof}
If $n\ge d$ then $\phi_{B,n}^0$ is an isomorphism. Indeed, one can show that the isomorphism of \cite[prop 1.3]{FF} equals $\phi_{B,n}^0$. Another proof of that fact, relying on classical invariant theory, is given in \cite[prop. 4.9.2]{Touze}. 

Moreover $B\mapsto H^*_{\mathcal{P}}(B)$ and $B\mapsto H^*(GL_{n,\mathbb{k}},B(\mathbb{k}^n,\mathbb{k}^n))$ are $\delta$-functors which vanish on the injectives (for $*>0$) and 
$\phi_{B,n}$ is a morphism of $\delta$-functors. Since $\phi_{B,n}^0$ is an isomorphism, $\phi_{B,n}^*$ must be an isomorphism \cite[p. 140]{Tohoku}.
\end{proof}

We now turn to examining the compatibility of the maps $\phi_{B,n}$ with the cup products. We first recall the definition of cross products and (external) cup products.
Let $(\mathfrak{A},\otimes)$ be an abelian category with enough injectives, equipped with a biexact monoidal product $\otimes$ which preserves the injectives. Then we may define a (associative, graded) ``cross product'' for $\Ext$-groups in $\mathfrak{A}$:
$$\times\,:\,\Ext^*_\mathfrak{A}(A_1,B_1)\otimes \Ext^*_\mathfrak{A}(A_2,B_2)\to \Ext^{*}_\mathfrak{A}(A_1\otimes A_2,B_1\otimes B_2)\;.$$
This cross product may be computed in two ways. First, using Yoneda
extensions. If $B_i\hookrightarrow E_i^\bullet$, $i=1,2$, are two extensions
representing classes $e_i\in\Ext^{k_i}_\mathfrak{A}(A_i,B_i)$, then $e_1\times
e_2$ is the class represented by the $(k_1+k_2)$-fold extension $B_1\otimes
B_2\hookrightarrow E^\bullet_1\otimes E_2^\bullet$. Second, using injective
coresolutions. For $i=1,2$, let $J_i^\bullet$ be an injective coresolution of
$B_i$. Then $J_1^\bullet\otimes J_2^\bullet$ is an injective coresolution of
$B_1\otimes B_2$. If $\alpha_i\in \hom(A_i,J_i^{k_i})$ are cocycles
representing classes $[\alpha_i]\in \Ext^{k_i}_\mathfrak{A}(A_i,B_i)$, then $[\alpha_1]\times [\alpha_2]$ is the class represented by the cocycle $\alpha_1\otimes \alpha_2\in \hom(A_1\otimes A_2,J_1^{k_1}\otimes J_2^{k_2})$.

Thus, if $\mathbb{k}$ is a field of positive characteristic then the category $(\mathcal{P}_{\mathbb{k}}(1,1),\otimes)$ of strict polynomial bifunctors and the category $(G_{\mathbb{k}}\text{-mod},\otimes)$ of rational modules over an algebraic group scheme $G_\mathbb{k}$ (both equipped with the usual tensor products) both have a cross product. Bifunctor cohomology is equipped with a (associative, graded) cup product: 
$$\begin{array}{cccc}
\cup:&\Ext^k_{\mathcal{P}_{d,\mathbb{k}}^d}(\Gamma^d gl,B)
\otimes \Ext^\ell_{\mathcal{P}_{e,\mathbb{k}}^e}(\Gamma^e gl,B')&
\to &\Ext^{k+\ell}_{\mathcal{P}_{d+e,\mathbb{k}}^{d+e}}(\Gamma^{d+e} gl,B\otimes B')\,,\\
&x\otimes y & \mapsto & \Delta_{d,e}^*(x\times y)
\end{array}
 $$
where $\Delta_{d,e}$ is the map induced by the diagonal $\Gamma^{d+e}\hookrightarrow \Gamma^d\otimes \Gamma^e$. Similarly, the rational cohomology of $G_{\mathbb{k}}$ is equipped with a cup product~:
$$\begin{array}{cccc}
\cup:&\Ext^k_{G_{\mathbb{k}}\text{-mod}}(\mathbb{k},M)\otimes \Ext^\ell_{G_{\mathbb{k}}\text{-mod}}(\mathbb{k},M')&\to &\Ext^{k+\ell}_{G_{\mathbb{k}}\text{-mod}}(\mathbb{k},M\otimes M')\,,\\
&x\otimes y & \mapsto & \Delta^*(x\times y)
\end{array}
 $$
where $\Delta$ is the isomorphism $\mathbb{k}\simeq \mathbb{k}\otimes \mathbb{k}$, $1\mapsto 1\otimes 1$. 
 
\begin{lemme}\label{lm-compat-phi-cup}
The natural map $\phi_{B,n}$ is compatible with the cup products:
$$\phi_{B\otimes C,n}(x\cup y)=\phi_{B,n}(x)\cup \phi_{C,n}(y)\;.$$
\end{lemme}
\begin{proof}
Evaluation on the pair $(\mathbb{k}^d,\mathbb{k}^d)$ is compatible with the cross products. Thus, the compatibility of $\phi_{B,n}$ with the cup products comes from the commutativity of the diagram:
$$\xymatrix{
\Gamma^{d+e}(gl(\mathbb{k}^n,\mathbb{k}^n))\ar@{->}[r]^-{\Delta_{d,e}}& 
\Gamma^d(gl(\mathbb{k}^n,\mathbb{k}^n))\otimes 
\Gamma^e(gl(\mathbb{k}^n,\mathbb{k}^n))\\
\mathbb{k}\ar@{->}[r]^{\Delta}\ar@{->}[u]&\mathbb{k}\otimes \mathbb{k}\ar@{->}[u]
}$$
where the vertical morphisms are defined using the $GL_{n,\mathbb{k}}$-equivariant maps $x\mapsto x\,(\Id_{\mathbb{k}^n})^{\otimes i}$, for $i=d,e$ and $d+e$.
\end{proof}

Putting all this information together we obtain the following strengthening of \cite[Th 1.5]{FF}:
\begin{Theorem}\label{thm-cohom-bif-computes-cohom-rat}
Let $\mathbb{k}$ be a field of positive characteristic and let $B\in \mathcal{P}_{d,\mathbb{k}}^d$ be a homogeneous bifunctor of bidegree $(d,d)$. For all $n\ge 1$, there are maps 
$$\phi_{B,n}:H^*_{\mathcal{P}}(B)\to H^*(GL_{n,\mathbb{k}},B(\mathbb{k}^n,\mathbb{k}^n))\;,$$
natural in $B$ and compatible with the cup products: $\phi_{B\otimes C,n}(x\cup y)=\phi_{B,n}(x)\cup \phi_{C,n}(y)$. Moreover, if $n\ge d$ then $\phi_{B,n}$ is an isomorphism.  
\end{Theorem}

\subsection{Proof of theorem \ref{thm-lifted-classes} assuming theorem \ref{thm-lifted-classes-cohom-bif}}

We now  prove that theorem \ref{thm-lifted-classes} is implied by the following bifunctor cohomology result:
\begin{Theorem}\label{thm-lifted-classes-cohom-bif}
Let $\mathbb{k}$ be a field of characteristic $p>0$. There are cohomology classes $c[d]\in H^{2d}_\mathcal{P,\mathbb{k}}(\Gamma^{d}(gl^{(1)}))$ such that~:
\begin{enumerate}
\item $c[1]\in H^2_\mathcal{P}(gl^{(1)})$ is non zero.
\item Let $d\ge 1$ and let $\Delta_{(1,\dots ,1)}:\Gamma^d (gl^{(1)})\to (gl^{(1)})^{\otimes d}$ be the map induced by the diagonal $\Gamma^d\to \otimes^d$. Then $\Delta_{(1,\dots ,1)\,*}c[d]=c[1]^{\cup d}$.
\end{enumerate}
\end{Theorem}

In order to prove theorem \ref{thm-lifted-classes}, we first remark that it is in fact a stable rational cohomology statement:
\begin{lemme}\label{lm-thm-lift-class-est-stable} Let $n_0$ be an integer greater or equal to the characteristic of $\mathbb{k}$.
Suppose that theorem \ref{thm-lifted-classes} is valid for $n=n_0$. Then theorem \ref{thm-lifted-classes} is valid for all $n$ such that $2\le n\le n_0$.
\end{lemme}
\begin{proof} The inclusion of $\mathbb{k}^n$ into the first $n$ coordinates of $\mathbb{k}^{n_0}$ and the projection $\mathbb{k}^{n_0}\to \mathbb{k}^n$ onto the first $n$ coordinates induce a map $\mathfrak{gl}_{n_0}\to \mathfrak{gl}_n$. Together with the inclusion $GL_{n,\mathbb{k}}\to GL_{n_0,\mathbb{k}}$, $M\mapsto [\text{\tiny$\begin{array}{cc}M & 0\\ 0 & 1 \end{array} $}]$, they induce `restriction' maps $H^*(GL_{n_0,\mathbb{k}}, \Gamma^m(\mathfrak{gl}_{n_0}^{(1)}))\to H^*(GL_{n,\mathbb{k}}, \Gamma^m(\mathfrak{gl}_{n}^{(1)})) $. These maps send the set of classes $c[m]\in H^{2m}(GL_{n_0,\mathbb{k}}, \Gamma^m(\mathfrak{gl}_{n_0}^{(1)}))$ to a set of classes $c'[m]\in H^{2m}(GL_{n,\mathbb{k}}, \Gamma^m(\mathfrak{gl}_{n}^{(1)}))$. By naturality of the restriction maps, the classes $c'[m]$ also satisfy condition 2 of theorem \ref{thm-lifted-classes}. To finish the proof, we have to check that $c'[1]$ is not null. The class $c[1]$ is not null, and by \cite{FS} $H^2(GL_{n_0,\mathbb{k}},\mathfrak{gl}_{n_0}^{(1)})$ is one dimensional,  generated by the Witt vector class. As remarked in \cite[remark 4.1]{VdK}, this implies that the restriction of $c[1]$, and hence of $c'[1]$, to an infinitesimal one parameter subgroup $G_{a1}$ is non trivial. Thus, $c'[1]$ is  non trivial. 
\end{proof}

\begin{proof}[proof of theorem \ref{thm-lifted-classes}]
Let's suppose that theorem  \ref{thm-lifted-classes-cohom-bif} is true. By lemma \ref{lm-thm-lift-class-est-stable}, it suffices to prove theorem \ref{thm-lifted-classes} for $n\ge p$. The maps $\phi_{\Gamma^m(gl^{(1)}),n}$ send the bifunctor cohomology classes $c[i]$ of theorem \ref{thm-lifted-classes-cohom-bif} to rational cohomology classes still denoted by $c[i]$. By naturality of the $\phi_{\Gamma^m(gl^{(1)}),n}$ and compatibility with the cup products, the rational cohomology classes $c[i]$ satisfy condition 2 of theorem \ref{thm-lifted-classes}. Since $n\ge p$, $\phi_{\Gamma^1(gl^{(1)}),n}$ is an isomorphism. Thus the rational cohomology class $c[1]$ is not null.
\end{proof}

\section{Complexes and $p$-complexes}\label{sec-p-complexes}

Homological algebra for $N$-complexes has been developed in \cite{KW,Kap}, and used for computations in quantum differential calculus. When $N=p$ is a prime, $p$-homological algebra is also the natural framework for some combinatorics of representation theory over fields of characteristic $p$ \cite{Troesch}. In this section, we recall the basic definitions and properties of $N$-complexes. When $N=p$, we prove a tensor product formula (proposition \ref{prop-comparaison-prod-tens}) which enable us to identify  cup products in $p$-coresolutions in section \ref{sec-proof-existence-lifted-classes}. To avoid confusion, we denote the $N$-complexes by the letters $C,D$ and the ordinary complexes by the letters $K,L$ in this section.

\subsection{Contractions of $N$-complexes}
Let $\mathfrak{A}$ be an additive category and let $N\ge 2$ be an integer. 
A $N$-complex  in $\mathfrak{A}$ is a graded object
$$C^{\bullet}=\bigoplus_{n\in\mathbb{N}} C^{n}$$
equipped with a $N$-differential, ie. a morphism $d$ of degree $1$ such that $d^N=0$.

For all integer $1\le s\le N-1$ we can `contract' a $N$-complex $C$ into an
ordinary complex $C_{[s]}$ by taking alternatively $d^s$ and $d^{N-s}$ as 
differentials~:
$$C_{[s]}:\;\; C^0\xrightarrow[]{d^s}C^s\xrightarrow[]{d^{N-s}}C^N\xrightarrow[]{d^s}C^{N+s}\xrightarrow[]{d^{N-s}}C^{2N}\to\dots $$
On the reverse way, we can build $N$-complexes out of ordinary complexes.
We define the $N$-complex $\widetilde{K}$ associated to an ordinary complex $K$
$$\widetilde{K}\;:\quad K^0\to \underbrace{K^1\xrightarrow[]{=}K^1\xrightarrow[]{=}\dots\xrightarrow[]{=}K^1}_{\text{$N-1$ terms}}\to K^2\to \underbrace{K^3\xrightarrow[]{=} K^3\xrightarrow[]{=}\dots}_{\text{$N-1$ terms}}\;.$$
We now specify the link between $N$-complexes and ordinary complexes.
\begin{lemme}\label{lm-2-res-p-res}
Let $K$ be an ordinary complex. For all $s\in [1,N-1]$ we have an equality $${(\widetilde{K})_{[s]}}=K\;.$$
Let $C$ be a $N$-complex. There is a morphism $\eta_C$ of $N$-complexes:
$$\eta_C: \widetilde{(C_{[1]})} \to C\;,$$
natural in $C$, such that if $N$ divides $i$ or $i-1$ then $\eta_C^i:\widetilde{(C_{[1]})}^i \to C^i$ is the identity.
\end{lemme}
\begin{proof}
The first claim follows directly from the definitions. To define $\eta_C$, we use the commutative diagram~:
$$\xymatrix @C=0.4cm{
\widetilde{(C_{[1]})}\,:& C^0\ar@{->}[r]^{\partial}\ar@{->}[d]^{=}& C^1\ar@{->}[r]^{=}\ar@{->}[d]^{=}&C^1\ar@{->}[r]^{=}\ar@{->}[d]^{\partial}&\dots\ar@{->}[r]^{=}&C^1\ar@{->}[r]^{\partial^{p-1}}\ar@{->}[d]^{\partial^{p-2}}& C^p\ar@{->}[r]^{\partial}\ar@{->}[d]^{=}&C^{p+1}\ar@{->}[d]^{=}\ar@{->}[r]^{=}&\dots \\
C\,:& C^0\ar@{->}[r]^{\partial}& C^1\ar@{->}[r]^{\partial}&C^2\ar@{->}[r]^{\partial}&\dots\ar@{->}[r]^{\partial}&C^{p-1}\ar@{->}[r]^{\partial}& C^p\ar@{->}[r]^{\partial}&C^{p+1}\ar@{->}[r]^{\partial}&\dots 
\;\;.}$$
\end{proof}

\subsection{Tensor product of $p$-complexes}
In this section, $p$ is a prime. We work in a $\mathbb{F}_p$-linear additive category $\mathfrak{A}$ equipped with a biexact monoidal product $\otimes$. 

Let $C$ and $D$ be two $p$-complexes. Since $p$ is a prime, we have 
$$(d_{C}\otimes 1 + 1 \otimes d_{D})^p = \sum_{i=0}^p {p\choose i}d_{C}^i\otimes d_{D}^{p-i}= d_{C}^p\otimes 1+1\otimes d_{D}^p=0\;.$$
Thus the differential $d_{C}\otimes 1 + 1 \otimes d_{D}$ (without sign!), 
makes the tensor product $C\otimes D$ into a $p$-complex. 

If $C,D$ are $p$-complexes, we have two ways of contracting and 
taking tensor products.
First, we may take the tensor product of the $p$-complexes $C$ and $D$, 
and take the ordinary complex $(C\otimes D)_{[1]}$ associated to this
$p$-complex for $s=1$. 
We may also consider the ordinary tensor product (with a sign!) 
$C_{[1]}\otimes D_{[1]}$ of the contracted complexes $C_{[1]}$ and $D_{[1]}$. 
In general, the complexes $(C\otimes D)_{[1]}$ and $C_{[1]}\otimes D_{[1]}$ 
are \emph{not equal}. For example, if $p=3$, the beginning of the complex 
$(C\otimes D)_{[1]}$ has the form:
$$ C^0\otimes D^0 \to 
\begin{array}{c}
C^0\otimes D^1\\
\oplus\; C^1\otimes D^0
\end{array}
\to 
\begin{array}{c}
C^0\otimes D^3\\
\oplus\; C^1\otimes D^2\\
\oplus\; C^2\otimes D^1\\
\oplus\; C^3\otimes D^0
\end{array}
\to 
\begin{array}{c}
C^0\otimes D^4\\
\oplus\; C^1\otimes D^3\\
\oplus\; C^2\otimes D^2\\
\oplus\; C^3\otimes D^1\\
\oplus\; C^4\otimes D^0
\end{array}
\to 
\dots
$$
while the beginning of the complex $C_{[1]}\otimes D_{[1]}$ has the form:
$$C^0\otimes D^0 \to 
\begin{array}{c}
C^0\otimes D^1\\
\oplus\; C^1\otimes D^0
\end{array}
\to 
\begin{array}{c}
C^0\otimes D^3\\
\oplus\; C^1\otimes D^1\\
\oplus\; C^3\otimes D^0
\end{array}
\to 
\begin{array}{c}
C^0\otimes D^4\\
\oplus\; C^1\otimes D^3\\
\oplus\; C^3\otimes D^1\\
\oplus\; C^4\otimes D^0
\end{array}
\to 
\dots\;.
$$
However, these two complexes have some similarities. For example we have~:
\begin{lemme}\label{lm-sous-obj}
Let $C$, $D$ be two $p$-complexes. For all nonnegative integers $k,\ell$, the object $C^{kp}\otimes D^{\ell p}$ appears once and only once in the complex $(C\otimes D)_{[1]}$ (resp. in the complex $C_{[1]}\otimes D_{[1]}$). Moreover, it appears in degree $2(k+\ell)$. 
\end{lemme}

\begin{Definition}Let $C$, $D$ be two $p$-complexes.
We define $p(C,D)^*$ as the graded object which equals 
$$p(C,D)^*:=\bigoplus_{k,\ell\,\ge 0}C^{kp}\otimes D^{\ell p}\;,$$
with $C^{kp}\otimes D^{\ell p}$ in degree $2(k+\ell)$. Lemma \ref{lm-sous-obj} yields inclusions of $p(C,D)^*$ into the graded objects $(C\otimes D)_{[1]}^*$ and $(C_{[1]}\otimes D_{[1]})^*$.
\end{Definition}

We now come to the main result of this section, which compares the complexes $(C\otimes D)_{[1]}$ and $C_{[1]}\otimes D_{[1]}$. We need this result in section \ref{sec-proof-existence-lifted-classes}.
\begin{Proposition}\label{prop-comparaison-prod-tens}
Let $C$, $D$ be two $p$-complexes. There is a morphism of ordinary complexes 
$$h_{C,D}: (C_{[1]}\otimes D_{[1]})\to (C\otimes D)_{[1]}$$
with the following properties~:
\begin{enumerate}
\item $h_{C,D}$ is natural with respect to the $p$-complexes $C,D$. 
\item $h_{C,D}^0$ and $h_{C,D}^1$ are identity maps.
\item there is a commutative diagram of graded objects~:
$$\xymatrix{
(C_{[1]}\otimes D_{[1]})^*\ar@{->}[rr]^-{h_{C,D}^*}&& (C\otimes D)_{[1]}^*\\
p(C,D)^*\ar@{^{(}->}[u]\ar@{=}[rr]&& p(C,D)^*\ar@{^{(}->}[u]\;.
}$$
\end{enumerate}
\end{Proposition}

We prove proposition \ref{prop-comparaison-prod-tens} in paragraph
\ref{subsec-proof} below. But before this, we indicate why proposition 
\ref{prop-comparaison-prod-tens} is useful.

\subsection{Homology of $p$-complexes and K\"unneth formulas}
We still work in a $\mathbb{F}_p$-linear category $\mathfrak{A}$
with a biexact monoidal product $\otimes$, but we now assume furthermore
that $\mathfrak{A}$ is \emph{abelian}. So, if $C$ is a $p$-complex in
$\mathfrak{A}$, the homology groups of the contracted complexes $C_{[s]}$,
$1\le s<p$, are well defined. 

\begin{Definition}A $p$-coresolution of $F\in\mathfrak{A}$ is a $p$-complex 
$C$ such that for all $s\in[1,N-1]$ the complex
$C_{[s]}$ is a coresolution of $F$.
A  $p$-acyclic complex is a $p$-coresolution of $0$.
\end{Definition}

The following weak K\"unneth formula is due to Troesch \cite[Th 2.3.1]{Troesch}. 

\begin{Proposition}\label{prop-tens-p-complexes}
Let $C$ be a $p$-coresolution of $F$ and let $D$ be a $p$-coresolution of $G$. The tensor product $\left(C\otimes D\right)$ is a $p$-coresolution of $F\otimes G$.
\end{Proposition}

If $C,D$ are $p$-complexes, we do not know any general K\"unneth formula 
relating the homology of $(C\otimes D)_{[s]}$ and
the tensor products of the homology of 
the contracted complexes $C_{[s]}$ and $D_{[s']}$. 
Our proposition \ref{prop-comparaison-prod-tens} is a substitute to the
K\"unneth formula: at least we know a natural chain map $h_{C,D}$ which allows
us to compare easily even degree cycles 
in $(C_{[1]}\otimes D_{[1]})$ and in $(C\otimes D)_{[1]}$. (If $C,D$ are 
$p$-coresolutions, propositions \ref{prop-tens-p-complexes} and
\ref{prop-comparaison-prod-tens}(2) show that $h_{C,D}$ is a quasi 
isomorphism).

\subsection{Proof of proposition
\ref{prop-comparaison-prod-tens}}\label{subsec-proof}
If the characteristic $p$ equals $2$, then the complexes
$C_{[1]}\otimes D_{[1]}$ and $(C\otimes D)_{[1]}$ are equal and there is nothing to prove. Therefore, we may suppose that $p$ is odd.

%
By lemma \ref{lm-2-res-p-res}, we have a morphism of complexes, natural in the $p$-complexes $C, D$~:
$$(\eta_C\otimes \eta_D)_{[1]} \,:\;(\widetilde{C_{[1]}}\otimes \widetilde{D_{[1]}})_{[1]}\to (C\otimes D)_{[1]}\;,$$
such that $(\eta_C\otimes \eta_D)_{[1]}^0$ and $(\eta_C\otimes \eta_D)_{[1]}^1$ are identity maps, and which fits into a commutative diagram of graded objects~:
$$\xymatrix{
(\widetilde{C_{[1]}}\otimes \widetilde{D_{[1]}})_{[1]}^*\ar@{->}[rr]^-{(\eta_C\otimes \eta_D)_{[1]}^*} && (C\otimes D)_{[1]}^*\\
p(\widetilde{C_{[1]}},\widetilde{D_{[1]}})^*
\ar@{^{(}->}[u]\ar@{=}[rr]
&& p(C,D)^*\ar@{^{(}->}[u]\;.
}$$
If $K,L$ are ordinary complexes, we define $p(K,L)^*$ as the graded subobject of $K^*\otimes L^*$ given by
$$p(K,L)^*:=
\bigoplus_{k,\ell\,\ge 0}K^{2k}\otimes L^{2\ell}\;.$$
Note that $p(K,L)^*= p(\widetilde{K},\widetilde{L})^*$. Thus, to prove proposition \ref{prop-comparaison-prod-tens}, it suffices to build a map of complexes $$H_{K,L}\;:\;  K\otimes L
\to (\widetilde{K}\otimes \widetilde{L})_{[1]}\;,$$
natural in the complexes $K, L$, which is the identity in degrees $0$ and $1$, and which fits into a commutative diagram~:
$$\xymatrix{
(K\otimes L)^*\ar@{->}[rr]^-{H_{K,L}^*}&& (\widetilde{K}\otimes \widetilde{L})_{[1]}^*\\
p(K,L)^*\ar@{^{(}->}[u]\ar@{=}[rr]&& p(\widetilde{K},\widetilde{L})^*\ar@{^{(}->}[u]\;.
}$$
Indeed, if such a map $H_{K,L}$ exists, then we may define $h_{C,D}$ as the composite $(\eta_C\otimes \eta_D)_{[1]}\circ H_{C_{[1]},D_{[1]}}$. In the remainder of this section, we give an explicit construction of the map $H_{K,L}$.

\subsubsection{Description of the complexes $K\otimes L$ and  $(\widetilde{K}\otimes \widetilde{L})_{[1]}$}

\begin{lemme}\label{lm-descr-cplx-ordin}
Let $K,L$ be two complexes. The objects of the complex $K\otimes L$ are given by the formulas~:
$$
(K\otimes L)^{2n} = T_{2n}\oplus T'_{2n}\;,\qquad (K\otimes L)^{2n+1} = T_{2n+1}\oplus T'_{2n+1}\;.$$
where the terms $T_{2n}$, $T'_{2n}$, $T_{2n+1}$ and $T'_{2n+1}$ are given by~:
\begin{align*}
&T_{2n}=\bigoplus_{k=0}^n  {K}^{2k}\otimes {L}^{2(n-k)}\;,&&&& T'_{2n}=\bigoplus_{k=1}^{n-1} {K}^{2k+1}\otimes {L}^{2(n-k)-1}\;,\\
&T_{2n+1}=\bigoplus_{k=0}^{n}  {K}^{2k}\otimes {L}^{2(n-k)+1}\;,&&&&
T'_{2n+1}=\bigoplus_{k=0}^{n}  {K}^{2k+1}\otimes {L}^{2(n-k)}\;.
\end{align*}
The differential $d$ of the complex $K\otimes L$ is given by the formula~: 
\begin{align*}
d(x) &= (d_{K}\otimes 1)(x)+ (1\otimes d_{L})(x)\quad\text{ if $x$ is in  $T_{2n}$ or $T_{2n+1}$,}\\
d(x') &= (d_{K}\otimes 1)(x')- (1\otimes d_{L})(x')\quad\text{ if $x'$ is in  $T'_{2n}$ or $T'_{2n+1}$.}
\end{align*}

\end{lemme}

\begin{proof}
This is just the classical definition of the usual tensor product of ordinary complexes.
\end{proof}

Now we examine the ordinary complex  $(\widetilde{K}\otimes \widetilde{L})_{[1]}$. If we draw the commutative diagram which defines the $p$-complex $\widetilde{K}\otimes \widetilde{L}$, we see the following pattern.
$$\text{\footnotesize $\begin{array}{|c|ccccc|c|ccccc|c|c}
&&&&&&&&&&&&&\\
\hline
04&&&14&&&24&&&34&&&44&\\
\hline
&&&&&&&&&&&&&\\
&&&&&&&&&&&&&\\
03&&&13&&&23&&&33&&&43&~\\
&&&&&&&&&&&&&\\
&&&&&&&&&&&&&\\
\hline
02&&&12&&&22&&&32&&&42&\\
\hline
&&&&&&&&&&&&&\\
&&&&&&&&&&&&&\\
01&&&11&&&21&&&31&&&41&~\\
&&&&&&&&&&&&&\\
&&&&&&&&&&&&&\\
\hline
00&&&10&&&20&&&30&&&40&\\
\hline
\end{array}$}$$
The big squares (like the one labelled ``$11$'') contains $(p-1)\times (p-1)$ objects and the small squares (like the one labelled ``$00$'') contains $1\times 1$ objects. Within a given square or rectangle labelled ``$ij$", the objects equal $K^i\otimes L^j$ and the $p$-differentials are identities. The differentials which go upwards from a rectangle or a square to another one equal $1\otimes d_{L}$, while the differentials which go towards the right equal $d_{K}\otimes 1$.

We keep the notations of lemma \ref{lm-descr-cplx-ordin}. The objects contained in the vertical rectangles are the objects of the $T_{2n+1}$, the objects contained in the horizontal rectangles are the objects of the $T'_{2n+1}$, the objects contained in the small squares are the objects of the $T_{2n}$ and the objects contained in the big squares are the objects of the $T'_{2n}$. Thus we have~:

\begin{lemme}\label{lm-descr-cplx-biz}
Let $K,L$ be two complexes. The objects of the (ordinary) complex $(\widetilde{K}\otimes \widetilde{L})_{[1]}$ are given by the formulas~:
\begin{align*}
&(\widetilde{K}\otimes \widetilde{L})_{[1]}^{2n} 
= T_{2n}\oplus {T'_{2n}}^{\oplus (p-1)}\;,\\
&(\widetilde{K}\otimes \widetilde{L})_{[1]}^{2n+1} 
= T_{2n+1}\oplus T'_{2n+1}\oplus {T'_{2n}}^{\oplus (p-2)}\;,
\end{align*}
with the terms $T_{2n}$, $T'_{2n}$, $T_{2n+1}$ and $T'_{2n+1}$ as defined in lemma \ref{lm-descr-cplx-ordin}.
\end{lemme}

In order to describe the differentials of the complex $(\widetilde{K}\otimes \widetilde{L})_{[1]}$ we need one more notation. We let $\delta_n$ be the ``signed diagonal morphism" 
$$\begin{array}{cccc}\delta_n :& {T'_{2n}}&\to& {T'_{2n}}^{\oplus (p-1)}\\
 &x&\mapsto &  (x,-x,x,-x,\dots,x,-x)
\end{array}\;.
$$
\begin{lemme}\label{lm-descr-cplx-biz-diff}
Let $K,L$ be two complexes.
The differential $\partial$ of the complex $(\widetilde{K}\otimes \widetilde{L})_{[1]}$ sends an element
$$(x,\delta_n(x'))\in T_{2n}\oplus {T'_{2n}}^{\oplus (p-1)}$$
of degree $2n$ to the element
$$\left(d(x-x')\;, 0\,\right)\;\in \left(T_{2n+1}\oplus T'_{2n+1}\right)\oplus {T'_{2n}}^{\oplus (p-2)}\;.$$
Here, $d$ is the differential of the complex $K\otimes L$ described in lemma \ref{lm-descr-cplx-ordin}. The differential $\partial$ sends an element 
$$(x,x',0)\in T_{2n+1}\oplus T'_{2n+1}\oplus {T'_{2n}}^{\oplus (p-2)}$$
of degree $2n+1$ to the element 
$$\left((1\otimes d_{L})(x)+(d_{K}\otimes 1)(x')\; ,\; -\delta_n(d_{K}\otimes 1)(x)+\delta_n(1\otimes d_{L})(x')\;\right)$$
in $T_{2n+2}\oplus {T'_{2n+2}}^{\oplus (p-1)}$.
\end{lemme}
\begin{proof}
Let $(x,\delta_n(x'))$ be an element of degree $2n$ in the complex $(\widetilde{K}\otimes \widetilde{L})_{[1]}$. This element may be represented as an element of degree $pn$ in the $p$-complex
$\widetilde{K}\otimes \widetilde{L}$:
$$\text{\footnotesize $\begin{array}{c|c|cccc|c|c}
&&&&&&&\\
\hline
&x&&&&&&\\
\hline
&&x'&&&&&\\
&&&-x'&&&&\\
~&&&&x'&&&~\\
&&&&&-x'&&\\
\hline
&&&&&&x&\\
\hline
&&&&&&&\\
\end{array}$}$$
The $p$-differentials in the big squares are identities.  As a result, the component of $\partial(x,\delta_n(x'))$ in the upper diagonals ${T'_{2n}}^{\oplus (p-2)}$ of the big squares is null and $\partial(x,\delta_n(x'))$ equals:
$$
(1\otimes d_{L})(x)+(d_{K}\otimes 1)(x)+(1\otimes d_{L})(x')+(d_{K}\otimes 1)(-x')\;.
$$
Since $x\in T_{2n}$ and $x'\in T'_{2n}$, lemma \ref{lm-descr-cplx-ordin} asserts that this element equals $d(x-x')$, where $d$ is the differential of the complex $K\otimes L$.

Let $(x,0,0)\in T_{2n+1}$ be a bihomogeneous element of degree $2n+1$ in the complex $(\widetilde{K}\otimes \widetilde{L})_{[1]}$. Then $\partial(x,0,0)$ is a sum of $p$ elements $(y_k)_{0\le k\le p-1}$ whose respective positions may be represented in the $p$-complex $\widetilde{K}\otimes \widetilde{L}$:
$$\text{\footnotesize $\begin{array}{c|c|cccc|c|c}
&&&&&&&\\
\hline
&y_0&&&&&&\\
\hline
&&y_1&&&&&\\
&&&\ddots&&&&\\
~&&&&\ddots&&&~\\
&x&&&&y_{p-1}&&\\
\hline
&&&&&&&\\
\hline
&&&&&&&\\
\end{array}$}$$
Since the $p$-differentials within the big square and within the vertical rectangle are identities, we compute that $y_0=(1\otimes d_{L})(x)$ and $y_{p-1}= (d_{K}\otimes 1)(x)$. Moreover, the formula $(d_{\widetilde{K}}\otimes 1+1\otimes d_{\widetilde{L}})^p(x)=0$ implies the following $p-2$ equalities: $$y_{p-1}+y_{p-2}=0,\, \dots ,\,y_2+y_1=0\;.$$
As a result, we have 
$$\partial(x,0,0)= (1\otimes d_{L})(x)-\delta_n(d_{K}\otimes 1)(x)\;.$$
The computation of $\partial(0,x',0)$ is similar.
\end{proof}

\subsubsection{Definition of $H_{K,L}$}
We let
\begin{align*}
&\begin{array}{cccc}
H^{2n}_{K,L}: & T_{2n}\oplus T'_{2n}&\to &T_{2n}\oplus {T'_{2n}}^{\oplus (p-1)}\\
& (x,x')&\mapsto &(x,\delta_n(-x'))
\end{array},
\end{align*}
and
\begin{align*}
&\begin{array}{cccc}
H^{2n+1}_{K,L}:&  T_{2n+1}\oplus T'_{2n+1}&\to & T_{2n+1}\oplus T'_{2n+1}\oplus {T'_{2n}}^{\oplus (p-2)}\\
&(x,x')&\mapsto & (x,x',0)\\
\end{array}.
\end{align*}
The graded map $H^{*}_{K,L} $ is natural with respect to the complexes $K$ and $L$. Moreover $H^{*}_{K,L} $ is the identity in degrees $0$ and $1$, and it fits into a commutative diagram of graded objects:
$$\xymatrix{
(K\otimes L)^*\ar@{->}[rr]^-{H_{K,L}^*}&& (\widetilde{K}\otimes \widetilde{L})_{[1]}^*\\
p(K,L)^*=T_{2n}\ar@{^{(}->}[u]\ar@{=}[rr]&& T_{2n}=p(\widetilde{K},\widetilde{L})^*\ar@{^{(}->}[u]\;.
}$$
The following lemma concludes the proof of proposition \ref{prop-comparaison-prod-tens}.
\begin{lemme}
$H^{*}_{K,L} $ induces a map of complexes 
$$H_{K,L}^\bullet\;:\;  (K\otimes L)^\bullet
\to {(\widetilde{K}\otimes \widetilde{L})_{[1]}}^\bullet\;.$$
\end{lemme}
\begin{proof} We have to show that $H_{K,L}$ commutes with the differentials.
Let $(x,x')\in T_{2n}\oplus T'_{2n}$ be an element of degree $2n$ of $K\otimes L$. Then 
\begin{align*}\partial(H^{2n}_{K,L}(x,x')) = \partial (x,\delta_n(-x') )=\left(d(x+x')\;, \;0\right)=H^{2n+1}_{K,L}(d(x,x'))\;.
\end{align*}
The first and the third equalities hold by definition of $H_{K,L}$ while the second equality follows from lemma \ref{lm-descr-cplx-biz-diff}.
Now let $(x,x')\in T_{2n+1}\oplus T'_{2n+1}$ be an element of degree $2n+1$ of $K\otimes L$. Then by lemma \ref{lm-descr-cplx-ordin}, $H^{2n+2}_{K,L}(d(x,x'))$ equals
\begin{align*}
 H^{2n+2}_{K,L}\left((d_{K}\otimes 1)(x)+(1\otimes d_{L})(x)+(d_{K}\otimes 1)(x')-(1\otimes d_{L})(x')\right)\,.
\end{align*}
The element $(d_{K}\otimes 1)(x)-(1\otimes d_{L})(x')$ lies in $T'_{2n+2}$ while the element $(1\otimes d_{L})(x)+(d_{K}\otimes 1)(x')$ lies in $ T_{2n+2}$. As a result, by definition of $H_{K,L}$,  the element $H^{2n+2}_{K,L}(d(x,x'))$ equals
\begin{align*}\left((1\otimes d_{L})(x)+(d_{K}\otimes 1)(x')\; ,\; -\delta_n(d_{K}\otimes 1)(x)+\delta_n(1\otimes d_{L})(x')\;\right)\,.
\end{align*}
But lemma \ref{lm-descr-cplx-biz-diff} tells us that this element equals $\partial (x,x',0)$, which by definition of  $H_{K,L}$ equals $\partial(H^{2n+1}_{K,L}(x,x'))$.
This concludes the proof that $H^\bullet_{K,L}$ is a chain map.
\end{proof}

\section{Building explicit coresolutions}\label{sec-explicit}
In this section, we develop methods to build explicit coresolutions of bifunctors, in order to compute their cohomology. We first notice in section \ref{sec-cohom-Fgl} that for bifunctors of the form $F(gl)$, it suffices to build injective resolutions of the functor $F$. Keeping this result in mind, we turn to building explicit injective coresolutions in the category of strict polynomial functors. We are interested in functors of the form $F(I^{(1)})$, that is, in functors obtained by precomposing a functor $F$ by the Frobenius twist $I^{(1)}$. In section \ref{sec-twist-comp-res}, we define the class of twist compatible coresolutions of strict polynomial functors. The name ``twist compatible'' is given after the following property: if $J_F$ is a twist compatible coresolution of $F$, we may build an explicit coresolution of the composite $F(I^{(1)})$ out of it. We don't know if all strict polynomial functors admit injective twist compatible resolutions, but we show in section \ref{sec-twist-comp-res-Gammagl} that the divided powers $\Gamma^d$ do.

\subsection{Bifunctor cohomology \emph{via} acyclic coresolutions}\label{sec-cohom-Fgl}

Let $\mathbb{k}$ be a field of positive characteristic, let $F\in\mathcal{P}_{d,\mathbb{k}}$ be a degree $d$ homogeneous strict polynomial functor over $\mathbb{k}$ and let $J$ be an injective coresolution of $F$ in $\mathcal{P}_{d,\mathbb{k}}$. We may precompose it by the bifunctor $gl$ to obtain a coresolution $J (gl)$ of the bifunctor $F(gl)$ in $\mathcal{P}_{d,\mathbb{k}}^d$. The objects of this coresolution are not injective, but the following lemma asserts that they are $H_\mathcal{P}^*$-acyclic.
\begin{lemme}\label{lm-acyclic-obj}
Let $J$ be an injective object of $\mathcal{P}_{d,\mathbb{k}}$. Then $H_\mathcal{P}^*(J(gl))=0$ if $*>0$. 
\end{lemme}
\begin{proof}
First, any injective $J\in\mathcal{P}_{d,\mathbb{k}}$ is a direct summand of a direct sum of standard injectives of the form $S^d (gl(\mathbb{k}^d,-))$ (use \cite[Th. 2.10]{FS} and duality \cite[Prop 2.6]{FS}). Thus, it suffices to show the $H_\mathcal{P}^*$-acyclicity of the bifunctor
$$S^d(gl(\mathbb{k}^d,gl(-_1,-_2)))\simeq S^d (gl(\mathbb{k}^d\otimes-_1,-_2))\;.$$
We use a theorem of Akin, Buchsbaum, Weyman \cite[Th. III.1.4]{ABW}. This theorem yields a filtration of the bifunctor $S^d(gl(\mathbb{k}^d\otimes-_1,-_2))$ whose associated graded object is the direct sum
$$\mathrm{Gr}\left(\,S^d (gl(\mathbb{k}^d\otimes-_1,-_2))\,\right) \simeq \bigoplus_{\lambda \text{ partition of weight $d$}} gl(W_\lambda(\mathbb{k}^d\otimes-_1),S_\lambda(-_2))\,. $$
Here, $S_\lambda$ is the Schur functor associated to the partition $\lambda$ and $W_\lambda=S_\lambda^\sharp$ is its dual \cite[Prop 2.6]{FS}. The bifunctors which appear in the direct sum are called `separable' and their cohomology is given \cite[Th. 1.5]{FF} in terms of extensions in $\mathcal{P}_{d,\mathbb{k}}$~:
$$H^*_\mathcal{P}(gl(W_\lambda(\mathbb{k}^d\otimes-_1),S_\lambda(-_2)) )\simeq \Ext^*_{\mathcal{P}_{d,\mathbb{k}}}(W_\lambda(\mathbb{k}^d\otimes-), S_\lambda(-))\;. $$
The extension groups which appear on the right are null if $*>0$ \cite[Fact 2.1]{Chalupnik}. As a result, the graded object associated to the filtration of $S^d(gl(\mathbb{k}^d\otimes-_1,-_2))$ is $H_\mathcal{P}^*$-acyclic. We deduce that $S^d(gl(\mathbb{k}^d\otimes-_1,-_2))$ is $H_\mathcal{P}^*$-acyclic.
\end{proof}

Since the coresolutions of the form $J(gl)$ have $H_\mathcal{P}^*$-acyclic objects, we may use them to compute the cohomology of the bifunctors of the form $F(gl)$ 
\cite[remarque 3 p. 148]{Tohoku}. More precisely:

\begin{Proposition}\label{prop-Fgl}
Let $F\in \mathcal{P}_{d,\mathbb{k}}$ be a strict polynomial functor of degree $d$. Let $J$ be an injective coresolution of $F$ in $\mathcal{P}_{d,\mathbb{k}}$ and let $K$ be an injective coresolution of $F(gl)$ in $\mathcal{P}_{d,\mathbb{k}}^d$.
Let $f:J (gl)\to K$ be a map of coresolutions over $\Id_{F(gl)}$. Then the morphism of complexes
$$\hom_{\mathcal{P}^d_{d,\mathbb{k}}}(\Gamma^d(gl),f)\,:\,\hom_{\mathcal{P}^d_{d,\mathbb{k}}}(\Gamma^d(gl),J(gl))\to \hom_{\mathcal{P}^d_{d,\mathbb{k}}}(\Gamma^d(gl),K)$$
induces an isomorphism in homology.
\end{Proposition}

We now specify how to compute cup products \emph{via} $H_\mathcal{P}^*$-acyclic coresolutions:

\begin{Proposition}\label{prop-cup-acyclic}
For $i=1,2$, let $J_i$ be an injective coresolution of $F_i\in \mathcal{P}_{d_i,\mathbb{k}}$, and
let $\alpha_i$ be a cocycle in the complex  $\hom(\Gamma^{d_i} (gl), J_i (gl))$, representing a class $[\alpha_i]\in H^*_{\mathcal{P}}(F_i(gl))$. Let $\Delta_{d_1,d_2}$ be the map obtained by evaluating the diagonal $\Gamma^{d_1+d_2}\hookrightarrow \Gamma^{d_1}\otimes \Gamma^{d_2}$ on the bifunctor $gl$.
Then 
$$\alpha_1\cup\alpha_2 := (\alpha_1\otimes \alpha_2)\circ \Delta_{d_1,d_2}\in \hom_{\mathcal{P}^{d_1+d_2}_{d_1+d_2,\mathbb{k}}}(\Gamma^{d_1+d_2} (gl), J_1 (gl)\otimes J_2 (gl))$$
is a cocycle representing the cohomology class $[\alpha_1]\cup[\alpha_2]$.
\end{Proposition}

\begin{proof} For $i=1,2$, let $K_i$ be an injective coresolution of the
bifunctor $F_i(gl)$ and let $f_i:J_i(gl)\to K_i$ be a morphism of
coresolutions over the identity map $F_i(gl)=F_i(gl)$.  The tensor product
$f_1\otimes f_2:J_1 (gl)\otimes J_2 (gl)\to K_1\otimes K_2$ is a map of
coresolutions over the identity map of $F_1(gl)\otimes F_2(gl)$. Moreover, it
sends $(\alpha_1\otimes \alpha_2)\circ \Delta_{d_1,d_2}$ to the cocycle
$((f_1\circ \alpha_1)\otimes (f_2\circ \alpha_2))\circ \Delta_{d_1,d_2}$. By
proposition \ref{prop-Fgl} and definition of the cup product the later cocycle  represents $[\alpha_1]\cup[\alpha_2]$ in the complex $\hom(\Gamma^{d_1+d_2}(gl), K_1\otimes K_2 )$. By proposition \ref{prop-Fgl} again, this means that $(\alpha_1\otimes \alpha_2)\circ \Delta_{d_1,d_2}$ represents $[\alpha_1]\cup[\alpha_2]$ in the complex $\hom(\Gamma^{d_1+d_2} (gl), J_1 (gl)\otimes J_2 (gl))$.
\end{proof}

\subsection{Twist compatible coresolutions}\label{sec-twist-comp-res}

In this paragraph we work in the category $\mathcal{P}_{\mathbb{k}}$ of 
strict polynomial functors over a field $\mathbb{k}$ of positive 
characteristic $p$. Thanks to the work of Troesch \cite{Troesch}, we know 
explicit injective $p$-coresolutions of the twisted injectives of 
$\mathcal{P}_{\mathbb{k}}$. However, these coresolutions are not natural. 
To make them natural, we need to restrict to a combinatorial subcategory 
$\mathcal{TP}_{\mathbb{k}}$ of $\mathcal{P}_{\mathbb{k}}$, called the twist 
compatible category. Finally we describe how to build an
$H^*_{\mathcal{P}}$-acyclic
coresolution of a twisted functor $F(gl^{(1)})$ from a twist compatible 
coresolution of $F$, that is from a coresolution of $F$ which lives in the 
twist compatible category. In order to define the twist compatible category, it is important not to use a categorical (ie: only `up to isomorphism') definition of the direct sum. For us, `the' direct sum of two strict polynomial functors $F$ and $G$ means the functor $F\oplus G$ which sends a vector space $V$ to the set of couples $(f,g)$ with $f\in F(V)$ and $g\in G(V)$. 

\subsubsection{The twist compatible category}

Let $\lambda=(\lambda_1,\dots,\lambda_n)$ be a $n$-tuple of positive integers. We denote by $S^\lambda$ the tensor product of symmetric powers: $S^\lambda:= \bigotimes_{i=1}^n S^{\lambda_i}$. Such strict polynomial functors are referred to as ``symmetric tensors''. They are injective objects of $\mathcal{P}_{\mathbb{k}}$. Let us denote by $p\lambda$ the $n$-tuple $p\lambda:=(p\lambda_1,\dots,p\lambda_n)$. The precomposition of $S^{\lambda}$ by the Frobenius twist $I^{(1)}$ yields a polynomial functor $S^{\lambda}(I^{(1)})$ and we have a monomorphism~:
$$\begin{array}{ccc}
S^{\lambda}(I^{(1)})&\hookrightarrow & S^{p\lambda}\\
\bigotimes_{i=1}^n(x_{i,1}^{(1)}\dots x_{i,\lambda_i}^{(1)} ) &\mapsto &  \bigotimes_{i=1}^n(x_{i,1}^{p}\dots x_{i,\lambda_i}^{p} )\;.
\end{array}
$$

\begin{Definition}\label{def-twist-compatible-map}Let $(\lambda^i)$ and $(\mu^j)$ be two finite families of tuples of positive integers and let $f:\bigoplus_i S^{\lambda^i} \to \bigoplus_j S^{\mu^j}$ be a morphism between two finite sums of symmetric tensors. We say that $f$ is \emph{twist compatible} if there exist a morphism $\overline{f}$ such that the following diagram commutes~:
$$\xymatrix{
\bigoplus_i S^{\lambda^i}(S^p)\ar@{->>}[d]\ar@{->}[rr]^-{f(S^p)}&& \bigoplus_j S^{\mu^j}(S^p)\ar@{->>}[d]\\
\bigoplus_i S^{p\lambda^i}\ar@{->}[rr]^-{\overline{f}}&& \bigoplus_j S^{p\mu^j}
}
$$
where the vertical epimorphisms are induced by the multiplications 
$$\begin{array}{ccc}
S^n(S^p) &\twoheadrightarrow & S^{np}\\
((x_{1,1}\dots x_{1,p})\dots(x_{n,1}\dots x_{n,p}))&\mapsto &(x_{1,1}\dots x_{1,p}\dots x_{n,1}\dots x_{n,p})\;.
\end{array}$$
\end{Definition}

\begin{lemme}\label{lm-lifting}
If $f$ is twist compatible, the morphism $\overline{f}$ is uniquely determined. Moreover, we have a commutative diagram~:
$$\xymatrix{
\bigoplus_i S^{\lambda^i}(I^{(1)})\ar@{^{(}->}[d]\ar@{->}[rr]^-{f(I^{(1)})}&& \bigoplus_j S^{\mu^j}(I^{(1)})\ar@{^{(}->}[d]\\
\bigoplus_i S^{\lambda^i}(S^p)\ar@{->>}[d]\ar@{->}[rr]^-{f(S^p)}&& \bigoplus_j S^{\mu^j}(S^p)\ar@{->>}[d]\\
\bigoplus_i S^{p\lambda^i}\ar@{->}[rr]^-{\overline{f}}&& \bigoplus_j S^{p\mu^j}
}$$
and the composite $S^{\lambda}(I^{(1)})\hookrightarrow S^{\lambda}(S^p)\twoheadrightarrow S^{p\lambda}$ equals the inclusion $S^{\lambda}(I^{(1)})\hookrightarrow S^{p\lambda}$.
\end{lemme}

\begin{lemme}\label{lm-comp-tw-compat}
Let $f$ and $g$ be twist compatible maps. The composite $f\circ g$ and the linear combinations $a f+ b g$ for all $a,b$ in the field $\mathbb{k}$ are twist compatible maps (resp. with $\overline{f\circ g}=\overline{f}\circ \overline{g}$ and $\overline{af+bg}=a\overline{f}+b \overline{g}$).
\end{lemme}

In general, one cannot say that the tensor product of two twist compatible maps is a twist compatible map. Indeed, tensors products of the form $(\bigoplus_i S^{\lambda^i})\otimes (\bigoplus_j S^{\mu^j})$ are not \emph{equal} to a direct sum of symmetric tensors, but only \emph{isomorphic} to it. So, the precise statement for tensor products is:
\begin{lemme}\label{lm-tens-tw-compat}
Let  $f:\bigoplus_i S^{\lambda^i}\to \bigoplus_k S^{\gamma^k}$ and $g: \bigoplus_j S^{\mu^j}\to \bigoplus_\ell S^{\nu^\ell}$ be two twist compatible maps, and let $\alpha:(\bigoplus_i S^{\lambda^i})\otimes (\bigoplus_j S^{\mu^j})\simeq \bigoplus_{i,j} S^{\lambda^i}\otimes S^{\mu^j}$ and $\beta:(\bigoplus_k S^{\gamma^k})\otimes (\bigoplus_\ell S^{\nu^\ell})\simeq \bigoplus_{k,\ell} S^{\gamma^k}\otimes S^{\nu^\ell}$ be the canonical isomorphisms. Then the composite $\beta\circ (f\otimes g)\circ \alpha^{-1}$ is a twist compatible map. 
\end{lemme}

\begin{lemme}\label{lm-puiss-sym-perm-twist-compat}
Let $i,j$ be two integers. The multiplication $m:S^i\otimes S^j\to S^{i+j}$ and the permutation $\tau:S^i\otimes S^j\to S^{j}\otimes S^i$ are twist compatible.
\end{lemme}

\begin{Remark}
The comultiplication  $\Delta:S^{i+j} \to S^i\otimes S^j$ is \emph{not} twist compatible in general.
\end{Remark}


We are now ready to define the twist compatible category $\mathcal{TP}_{\mathbb{k}}$. We want this category to contain the direct sums of symmetric tensors and the twist compatible maps. We also want it to be stable under tensor products, so we have to introduce the ``iterated symmetric tensors''. A $0$-iterated symmetric tensor is just a symmetric tensor, that is, a functor of the form $S^\lambda$, where $\lambda$ is a tuple of positive integers. For $n\ge 1$, a $n$-iterated symmetric tensor is a functor $F$ of the form
$F:= \bigotimes_{i=1}^k\bigoplus_{j=1}^\ell S_{i,j}$
where the $S_{i,j}$ are $(n-1)$-iterated symmetric tensors. 
If $F$ is a $n$-iterated symmetric tensor then we have a canonical isomorphism $$F:=\bigotimes_{i=1}^k\bigoplus_{j=1}^\ell S_{i,j}\simeq \bigoplus_{(j_1,\dots,j_k)\in \mathbb{N}^\ell} \bigotimes_{m=1}^k S_{m,j_m}  $$
between $F$ and a direct sum of $(n-1)$-iterated symmetric tensors. Composing such isomorphisms, we may associate to each iterated symmetric tensor $F$ an isomorphism $\xi_F:F\to F_0$ from $F$ to a direct sum of symmetric tensors $F_0$.

\begin{Definition}
The \emph{twist compatible category} $\mathcal{TP}_{\mathbb{k}}$ is the subcategory  of $\mathcal{P}_{\mathbb{k}}$ whose objects are the iterated symmetric tensors and whose morphisms are the maps $f:F\to G$ such that the composite $f_0:= \xi_G\circ f\circ {\xi_F}^{-1}$ is twist compatible.
\end{Definition}

\begin{lemme}\label{lm-twcat-stable}
The twist compatible category  $\mathcal{TP}_{\mathbb{k}}$ is an additive subcategory of $\mathcal{P}_{\mathbb{k}}$, stable under tensor products. 
\end{lemme}
\begin{proof}
By lemma \ref{lm-comp-tw-compat}, $\mathcal{TP}_{\mathbb{k}}$ is an additive subcategory of $\mathcal{P}_{\mathbb{k}}$. By Lemma \ref{lm-tens-tw-compat}, it is stable under tensor products. 
\end{proof}

\subsubsection{Natural injective $p$-coresolutions}

Thanks to the work of Troesch \cite{Troesch}, we know some explicit injective $p$-coresolutions of the twisted symmetric powers. These coresolutions generalize in characteristic $p$ odd the coresolutions previously known in characteristic $p=2$ \cite{FLS}. Let's denote by $I^{\oplus p}$ the $p$-times iterated direct sum of the identity functor $I$ (ie: $I^{\oplus p}$ sends a vector space $V$ to $V^{\oplus p}$). Then we have \cite[Th. 1 and prop. 3.2.1]{Troesch}:
\begin{Theorem}\label{thm-Troesch} The functor
$$B_n:=S^n(I^{\oplus p})\simeq\bigoplus_{i_0+\dots+i_{p-1}=n}S^{i_0}\otimes\dots\otimes S^{i_{p-1}}$$
is equipped with a $p$-differential $d$ such that~:
\begin{itemize}
\item[(1)] The cohomological degree of $S^{i_0}\otimes\dots\otimes S^{i_{p-1}}$ is $$0.i_0+1.i_1+\dots+(p-1)i_{p-1}$$
and the $p$-differential increases the cohomological degree by one.
\item[(2)] If $n$ is a multiple of $p$ then $\left(B_n^\bullet,d\right)$ is a $p$-coresolution of the twisted symmetric power $S^{n/p}(I^{(1)})$. Otherwise $\left(B_n^\bullet,d\right)$ is $p$-acyclic.
\item[(3)] The canonical isomorphism~:
$$B_*(V\oplus W)= S^*((V\oplus W)^{\oplus p})\simeq S^*(V^{\oplus p})\otimes S^*(W^{\oplus p})=B_*(V)\otimes B_*(W)$$
is an isomorphism of $p$-complexes.
\end{itemize}
\end{Theorem}

\begin{Corollary}\label{cor-p-res}
Let $\mu=(\mu_1,\dots,\mu_k)$ be a tuple of positive integers. Then the tensor product $B_{p\mu} := \bigotimes_{i=1}^k B_{p\mu_i}$ is a $p$-complex such that
\begin{itemize}
\item[1.] $B_{p\mu}$ is a $p$-coresolution of $S^{\mu}(I^{(1)})$, which equals $S^{p\mu}$ in degree $0$.
\item[2.] Let $f:S^\mu\to S^\nu$ be a map between symmetric tensors and let $\widetilde{f}: S^{p\mu}\to S^{p\nu}$ be a map which fits into the commutative diagram~:
$$\xymatrix{
S^{p\mu}\ar@{->}[rr]^-{\widetilde{f}}&& S^{p\nu}\\
S^{\mu}(I^{(1)})\ar@{^{(}->}[u]\ar@{->}[rr]^-{{f}(I^{(1)})}&& S^{\nu}(I^{(1)})\ar@{^{(}->}[u]\;.
} 
$$
Then $\widetilde{f}(I^{\oplus p}):B_{p\mu} = S^{p\mu}(I^{\oplus p})\to  S^{p\nu}(I^{\oplus p})=B_{p\nu}$ is a map of $p$-complexes over $f(I^{(1)})$, which agrees with $\widetilde{f}$ in degree $0$.
\end{itemize}
\end{Corollary}
\begin{proof}
The first point follows from theorem \ref{thm-Troesch}(2) and proposition \ref{prop-tens-p-complexes}. Let's prove the second point. If $\psi$ is a multiplication $\psi:S^k\otimes S^\ell\to S^{k+\ell}$, a comultiplication $\psi:S^{k+\ell}\to S^k\otimes S^\ell$ or a permutation $\psi:S^k\otimes S^\ell\simeq S^\ell\otimes S^k$, then the precomposition of $\psi$ by $I^{\oplus p}$ induces a map of $p$-complexes thanks to theorem \ref{thm-Troesch}(3). Any map between symmetric tensors is built out of multiplications, permutations and comultiplications \cite[p. 779]{Chalupnik}. As a result, if $\psi: S^\gamma\to S^\delta$ is a map between symmetric tensors, the precomposition of $\psi$ by $I^{\oplus p}$ induces a morphism of $p$-complexes $\psi(I^{\oplus p}): S^\gamma(I^{\oplus p})\to S^\delta(I^{\oplus p})$. Moreover, thanks to theorem \ref{thm-Troesch}(1), this map of $p$-complexes equals $\psi: S^\gamma\to S^\delta$ in cohomological degree $0$. Now suppose that $\gamma=p\lambda,\delta=p\mu$ and $\psi=\widetilde{f}$ lifts $f(I^{(1)}):S^{\mu}(I^{(1)})\to S^{\nu}(I^{(1)})$. Then the map of $p$-complexes $\widetilde{f}(I^{\oplus p}): B_{p\mu}\to B_{p\nu}$ is a lifting of $f(I^{(1)})$.
\end{proof}

By corollary \ref{cor-p-res}, we may associate to each symmetric tensor $S^\lambda$ an injective $p$-coresolution $S^{p\lambda}(I^{\oplus p})$ of $S^\lambda(I^{(1)})$ and to each map $f:S^\lambda\to S^\mu$ a map of $p$-coresolutions $\widetilde{f}(I^{\oplus p})$ over $f(I^{(1)}):S^\lambda(I^{(1)})\to S^\mu(I^{(1)})$. Unfortunately, there is no natural choice of $\widetilde{f}$ in general. Indeed, let $p=2$ and let $\tau$ be the transposition of $\Si_2$ which acts on $\otimes^2$ by permutation of the factors of the tensor product. Then the map $(1+\tau):\otimes^2\to \otimes^2$ equals the composite 
$\otimes^2\xrightarrow[]{\mathrm{mult}}S^2\xrightarrow[]{\mathrm{diag}}\otimes^2$. As a result we have at least two different liftings for the map $(1+\tau)$, namely $(1+\tau)(S^2):S^2\otimes S^2\to S^2\otimes S^2$ and the composite 
$S^2\otimes S^2\xrightarrow[]{\mathrm{mult}}S^4\xrightarrow[]{\mathrm{diag}}S^2\otimes S^2$. To obtain a natural choice of the lifting $\widetilde{f}$ we restrict to the twist compatible category $\mathcal{TP}_{\mathbb{k}}$.
\begin{Proposition}
There is an additive functor $$T:\mathcal{TP}_{\mathbb{k}}\to
\text{$p$-}Ch^{\ge 0}(\mathcal{P}_{\mathbb{k}})$$
from the twist compatible category to the category of (nonnegative) $p$-cochain complexes, which sends an object $F$ to an injective $p$-coresolution of $F(I^{(1)})$ and a twist compatible map $f$ to a map of $p$-coresolutions over $f(I^{(1)})$.
Moreover, there is a natural isomorphism $T(F)\otimes T(G)\simeq T(F\otimes G)$ over the identity map $ F(I^{(1)})\otimes G(I^{(1)})=(F\otimes G)(I^{(1)})$.
\end{Proposition}
\begin{proof} We first define $T$ on direct sums of symmetric tensors:
we send a direct sum of symmetric tensors $\bigoplus_i S^{\lambda_i}$ to the $p$-coresolution $\bigoplus_i B_{p\lambda_i}$ and a twist compatible map $f$ to the map $\overline{f}(I^{\oplus p})$. Thanks to lemma \ref{lm-lifting}, this map $\overline{f}$ is unique. Lemma \ref{lm-comp-tw-compat} shows the functoriality of $T$ as well as its additivity.

Now we extend $T$ to the whole category $\mathcal{TP}_{\mathbb{k}}$. If $F$ is a $n$-iterated symmetric tensor, we have a well defined isomorphism $\xi_F:F\to F_0$ onto a direct sum of symmetric tensors $F_0$ and we set $T(F):= T(F_0)$. If $f:F\to G$ is a map between iterated symmetric tensors, then $f_0:= \xi_G\circ f\circ {\xi_F}^{-1}$ is a twist compatible map and we set $T(f):= T(f_0)$.

If $F$ and $G$ are iterated symmetric tensors, then $\xi_{F\otimes G}$ equals the composite $\xi_{F_0\otimes G_0}\circ (\xi_F\otimes \xi_G)$. By definition, $T(F)\otimes T(G)$ equals $T(F_0)\otimes T(G_0)$, $T(F\otimes G)$ equals $T((F\otimes G)_0)$ and $\xi_{F_0\otimes G_0}(I^{(1)})$ induces an isomorphism $T(F_0)\otimes T(G_0)\to T((F\otimes G)_0)$. Thus, we have a natural isomorphism $T(F)\otimes T(G)\to T(F\otimes G)$ over $(\xi_{F\otimes G}^{-1} \circ \xi_{F_0\otimes G_0}\circ (\xi_F\otimes \xi_G))(I^{(1)})$ that is, over the identity map $ F(I^{(1)})\otimes G(I^{(1)})=(F\otimes G)(I^{(1)})$. 
\end{proof}

As a particular case of the construction of the functor $T$, we record the description of the $\Si_d$-module $T(\otimes^d)$ for further use:
\begin{lemme}\label{lm-descrip-sig-mod-Ttens}
Let $d$ be a positive integer. Then $T(\otimes^d)=T(S^1)^{\otimes d}$. Moreover, if $\sigma\in\Si_d$ acts on $\otimes^d$ by permuting the factors of the tensor product, then the map $T(\sigma):T(\otimes^d)\to  T(\otimes^d)$ equals the map $\sigma(T(S^1)):T(S^1)^{\otimes d}\to T(S^1)^{\otimes d}$.
\end{lemme}
\begin{proof}
The formula $T(\otimes^d)=T(S^1)^{\otimes d}$ follows from the definition of $T$ on symmetric tensors. If $\sigma\in\Si_d$, then $\sigma$ is twist compatible and $\overline{\sigma}=\sigma(S^p)$. Thus, $T(\sigma)=\overline{\sigma}(I^{\oplus p})=\sigma(T(S^1))$.
\end{proof}

\subsubsection{Twist compatible coresolutions}

If $\mathfrak{A}$ is an additive category, we denote by 
$Ch^{\ge 0}(\mathfrak{A})$ (resp. $\text{bi-}Ch^{\ge 0}(\mathfrak{A})$) the
category of nonnegative cochain complexes (resp. bicomplexes) in $\mathfrak{A}$.
\begin{Definition}
A twist compatible complex is an object of $Ch^{\ge
0}(\mathcal{TP}_\mathbb{k})$. In particular, a coresolution $J$ of $F\in
\mathcal{P}_\mathbb{k}$ is twist compatible if all its objects $J^k$ and all its 
differentials  $\partial^k:J^k\to J^{k+1}$ belong to the twist compatible 
category $\mathcal{TP}_{\mathbb{k}}$.
\end{Definition}

Let $J$ be a twist compatible complex. 
We first apply objectwise the functor $T$ to $J$. 
We obtain a commutative diagram
$$ T(J^0)\to T(J^1)\to \dots \to T(J^k)\to T(J^{k+1})\to \dots $$
The rows of this diagram are ordinary complexes since 
$T(\partial)\circ T(\partial)=T(\partial\circ\partial)=0$, while for all $k$,
 the $k$-th column $T(J^k)$ is an injective $p$-coresolution of 
 $J^k(I^{(1)})$. Second, we apply the functor $-_{[1]}$ columnwise. 
 Thus we obtain a bicomplex
$$ T(J^0)_{[1]}\to T(J^1)_{[1]}\to \dots \to T(J^k)_{[1]}\to T(J^{k+1})_{[1]}\to \dots  $$
\begin{lemme}\label{lm-twcomptwres}
If $J$ is a twist compatible coresolution of $F\in\mathcal{P}_{\mathbb{k}}$,
then the totalization of the bicomplex $T(J)_{[1]}$ is an injective
coresolution of $F(I^{(1)})$.
\end{lemme}
\begin{proof}
We compute the homology of the bicomplex $T(J)_{[1]}$ in two steps: first
along the columns and second along the rows. We obtain $F(I^{(1)})$
concentrated in bidegree $(0,0)$. Whence the result.
\end{proof}

Finally we may precompose all the objects of the bicomplex $T(J)_{[1]}$ by the bifunctor $gl$.
\begin{Definition}
We denote by $A$ the functor
$$ A\,:\, Ch^{\ge 0}(\mathcal{TP}_\mathbb{k}) \to \text{bi-}Ch^{\ge
0}(\mathcal{P}_\mathbb{k}(1,1))$$
which sends a twist compatible complex $J$ to the (first quadrant)
bicomplex $A(J)= T(J)_{[1]}(gl)$.
\end{Definition}

The letter `$A$' stands for `acyclic' (coresolution). Indeed lemmas
\ref{lm-acyclic-obj}  and
\ref{lm-twcomptwres}, and proposition \ref{prop-Fgl} yield: 

\begin{Proposition}\label{prop-twist-compat-res-donne-bicompl}
Let $J^\bullet$ be a twist compatible coresolution of
$F\in\mathcal{P}_\mathbb{k}$. 
Then the total complex associated to the bicomplex $A(J)$ 
is an $H^*_\mathcal{P}$-acyclic coresolution of $F(gl^{(1)})$. In particular,
the homology of the complex $H^0_\mathcal{P}(\mathrm{Tot}(A(J)))$ computes the
cohomology of the bifunctor $F(gl^{(1)})$.
\end{Proposition}

\subsection{A twist compatible coresolution of $\Gamma^n$}\label{sec-twist-comp-res-Gammagl}

Let $\mathbb{k}$ be a field of positive characteristic. The reduced bar  construction yields a functor~:
$$\overline{B} : \{\text{ CDGA-algebras }\} \to \{\text{ CDGA-algebras }\} $$ 
from the category of Commutative Differential Graded Augmented algebras over $\mathbb{k}$ to itself \cite[Chap. X]{ML}. 

Let's recall some classical examples associated with the reduced bar construction. Let $V$ be a finite dimensional $\mathbb{k}$-vector space. Let  $S(V)$, $\Lambda(V)$ and $\Gamma(V)$ be the  symmetric, exterior and divided powers algebras over $V$. If we define the degree of an element $v\in V$ to be respectively $0$, $1$, and $2$ and the differential to be  zero, then these algebras become CDGA-algebras.
Moreover, the injective morphisms~: 
\begin{align*}&\Lambda^n(V)\hookrightarrow V^{\otimes n} = S^1(V)^{\otimes n}\subset \overline{B}_n(S^*(V))\\
&\Gamma^n(V)\hookrightarrow V^{\otimes n} = \Lambda^1(V)^{\otimes n}\subset \overline{B}_{2n}(\Lambda^*(V))
\end{align*}
define maps of CDGA-algebras $\Lambda^*(V)\hookrightarrow \overline{B}(S^*(V))$ and $\Gamma^*(V)\hookrightarrow \overline{B}(\Lambda^*(V))$. The following is well known (see for example \cite[section 5.4]{Real}):
\begin{lemme}
The maps $\Lambda^*(V)\hookrightarrow \overline{B}(S^*(V))$, and $\Gamma^*(V)\hookrightarrow \overline{B}(\Lambda^*(V))$, as well as the composite:
$$\Gamma^*(V)\hookrightarrow \overline{B}(\Lambda^*(V))\hookrightarrow \overline{B}(\overline{B}(S^*(V))$$
are quasi isomorphisms.
\end{lemme}
\begin{proof}
To prove that the first two maps are quasi isomorphisms, use \cite[4 p. 02]{CartanSem} to reduce to a one dimensional vector space and then compute. To prove that the composite is also a quasi isomorphism, use that $\overline{B}$ preserves quasi isomorphisms \cite[X Th. 11.2]{ML}.  
\end{proof}

If $X=S,\Lambda,\Gamma$, then the multiplications $X^d(V)\otimes X^e(V)\to X^{d+e}(V)$, the diagonals $X^{d+e}(V)\to X^d(V)\otimes X^e(V)$ and the permutations $X^d(V)\otimes X^e(V)\simeq X^e(V)\otimes X^d(V)$ involved in the definition of the Hopf algebra structure on $X^*(V)$ are actually maps of strict polynomial functors. As a consequence, we may interpret the reduced bar construction and the quasi-isomorphism $\Gamma^*(V)\hookrightarrow  \overline{B}(\overline{B}(S^*(V))$ in the category $\mathcal{P}_{\mathbb{k}}$ of strict polynomial functors over $\mathbb{k}$. This category splits as a direct sum  $\mathcal{P}_{\mathbb{k}}=\bigoplus_{d\ge 0}\mathcal{P}_{d,\mathbb{k}}$, where $\mathcal{P}_{d,\mathbb{k}}$ is the subcategory of homogeneous strict polynomial functors of polynomial degree $d$.
We want to examine more carefully the homogeneous part of polynomial degree $d$ of the complex $\overline{B}(\overline{B}(S^*(-))$.

In order to do this, we first recall the construction \cite[X 10]{ML} of the chain complex $\overline{B}_\bullet(A)$. Let $I_\bullet$ be the kernel of the augmentation $\epsilon:A_\bullet\to \mathbb{k}$ and let $(sI)_\bullet$ be the suspension of $I_\bullet$. That is, $(sI)_\bullet$ is the complex defined by $(sI)_n=I_{n-1}$ and $(sd)_n=-d_{n-1}$. We denote by $[a_1|\dots|a_n]$ an element $a_1\otimes\dots\otimes a_n$ of the complex $(sI)^{\otimes n}$. For each $n\ge 1$, we define a chain map $d_E:s^{-n}((sI)^{\otimes n})\to s^{-(n-1)}((sI)^{\otimes n-1})$ by the formulas:
$$d_E([a_1|\dots|a_n])=\sum_{i=1}^{n-1}(-1)^{\epsilon_i}[a_1|\dots|a_{i-1}|a_i a_{i+1}|\dots |a_n]\,,\quad d_E([a_1])= 0\,, $$
with $\epsilon_i=\deg [a_1|\dots|a_i]$ and the complex $(sI)^{\otimes 0}$ equals $\mathbb{k}$ concentrated in degree $0$.
Since $d_E\circ d_E=0$ we have a (first quadrant) bicomplex: 
$$\mathbb{k}={(sI)}^{\otimes 0}\xleftarrow[]{d_E}\dots \xleftarrow[]{d_E}s^{-(n-1)}((sI)^{\otimes n-1})\xleftarrow[]{d_E} s^{-n}((sI)^{\otimes n})\xleftarrow[]{d_E}\dots  \;. $$
The reduced bar construction $\overline{B}_\bullet(A)$ is the total complex associated to this bicomplex. We are now ready to prove:

\begin{lemme}\label{lm-details-doublebar} The component of homogeneous polynomial degree $1$ of the complex $\overline{B}_{\bullet}(\overline{B}(S^*(-))$ equals $S^1$ concentrated in degree $2$. For $d\ge 2$, the component of homogeneous polynomial degree $d$ of the complex $\overline{B}_{\bullet}(\overline{B}(S^*(-))$ is:
$$\dots\leftarrow\underbrace{\bigoplus_{k=0}^{d-2}(\otimes^k)\otimes(\otimes^2)\otimes(\otimes^{d-k-2})}_{\text{degree $2d-1$}}  \xleftarrow[]{\prod(1-\tau_k)}  \underbrace{\otimes^d}_{\text{degree $2d$}} \leftarrow \underbrace{0\leftarrow 0\leftarrow \dots}_{\text{degrees $n>2d$}}\;,$$
where $\tau_k\in \Si_d$ is the transposition which exchanges $k+1$ and $k+2$.
\end{lemme}
\begin{proof}
We have $\overline{B}(S^*)=\bigoplus S^{n_1}\otimes\dots\otimes S^{n_k}$, where the sum is taken over all $k$-tuples $(n_1,\dots,n_k)$ of positive integers, for all $k\ge 0$ (with the convention that the $0$-tuple $(\;)$ corresponds to the constant term $\mathbb{k}$). An element $[s_1|\dots|s_k]\in S^{n_1}\otimes\dots\otimes S^{n_k}$ has degree $k$ and polynomial degree $\sum n_i$. 

Let $I$ be the kernel of the augmentation $\epsilon : \overline{B}(S^*)\to \mathbb{k}$ and let $_i I$ be the component of homogeneous polynomial degree $i$ of the complex $I$. Then $_0 I$ is null, $_1 I$ equals $S^1$ concentrated in degree $1$, and $_2 I$ equals 
$$\dots\leftarrow0\leftarrow S^2\leftarrow S^1\otimes S^1\leftarrow 0\leftarrow\dots\;,$$ with $S^1\otimes S^1$ placed in degree $2$ and $S^2$ placed in degree $1$. 

We first analyze the homogeneous component of polynomial degree $d=1$ of $\overline{B}_{\bullet}(\overline{B}(S^*(-))$. We recall that
if $F$ and $G$ are strict polynomial functors homogeneous of polynomial degrees $f,g$ then the tensor product $F\otimes G$ is homogeneous of polynomial degree $f+g$. Since $_0 I$ is null, this implies that the polynomial  degree $1$ homogeneous part of $s^{-n}((sI)^{\otimes n})$ is null, except if $n=1$. In the latter case, it equals $_1 I$. Thus the component of polynomial degree $1$ of the bicomplex defining $\overline{B}_{\bullet}(\overline{B}(S^*(-))$ equals $S^1$ concentrated in bidegree $(1,1)$ and we are done.

We perform a similar analysis in polynomial degree $d\ge 2$. Since $_0 I$ is null, the component of polynomial degree $d$ of $s^{-n}((sI)^{\otimes n})$ is null if $n>d$. If $n=d$ the component of homogeneous polynomial degree $d$ of $s^{-n}((sI)^{\otimes n})$ equals $s^{-d}((s\,_1I)^{\otimes d})$, that is, it equals $(S^1)^{\otimes d}$ placed in degree $d$. Finally, if $n=d-1$ the homogeneous polynomial degree $d$ part of $s^{-n}((sI)^{\otimes n})$ equals: $s^{-(d-1)} (\bigoplus_{k=0}^{d-2} (s\,_1I)^{\otimes k}\otimes s\,_2I\otimes (s\,_1I)^{\otimes d-k-2})$. As a consequence, the component of homogeneous polynomial degree $d$ of the bicomplex
$$\dots \xleftarrow[]{d_E}s^{-(d-1)}((sI)^{\otimes d-1})\xleftarrow[]{d_E} s^{-d}((sI)^{\otimes d})\xleftarrow[]{d_E}s^{-(d+1)}((sI)^{\otimes d+1})\xleftarrow[]{d_E}\dots  \;. $$
may be written as~:
$$\xymatrix{
\dots &\ar@{->}[l] \bigoplus_{k=0}^{d-2} (S^1)^{\otimes k}\otimes (S^1)^{\otimes 2}\otimes (S^1)^{\otimes d-k-2}\ar@{->}[d]&\ar@{->}[l]_-{d_E} (S^1)^{\otimes d}\ar@{->}[d]&0 \ar@{->}[l]\ar@{->}[d]&\dots \ar@{->}[l]\\
\dots &\ar@{->}[l] \bigoplus_{k=0}^{d-2} (S^1)^{\otimes k}\otimes S^2\otimes (S^1)^{\otimes d-k-2}&\ar@{->}[l] 0&0 \ar@{->}[l]&\dots \ar@{->}[l]
}$$
with $(S^1)^{\otimes d}$ placed in bidegree $(d,d)$. We now turn to showing that the map 
$d_E$ has the appropriate form. By definition, $d_E$ sends $[s_1|\dots |s_d]\in (S^1)^{\otimes d}$ to $\sum_{k=0}^{d-2}(-1)^{0}[s_1|\dots|s_k*s_{k+1}|\dots s_n]$, where $*$ denotes the multiplication in $\overline{B}(S^*)$. This multiplication is the ``shuffle product''
\cite[formula (12.4) p.313]{ML}. By definition, this shuffle product sends an element $[s_{k}]\otimes [s_{k+1}]\in \overline{B}_1(S^*)\otimes \overline{B}_1(S^*)$ to the sum $[s_{k}|s_{k+1}]-  [s_{k+1}|s_{k}]$. This concludes the proof.
\end{proof}

Let $J_d^\bullet$ be the homogeneous part of polynomial degree $d$ 
of the complex $\overline{B}_{2d-\bullet}(\overline{B}(S^*(-))$. 
We now state the main result of this section, compare \cite{Totaro}: 
\begin{Proposition}\label{prop-res-Gamma}Let $\mathbb{k}$ be a field of characteristic $p>0$. There is a  family $\left(J_d^\bullet\right)_{d\ge 1}$ of \emph{twist compatible} coresolutions of the divided powers $\Gamma^d$ such that $J_1$ equals $S^1$ concentrated in degree $0$, and for all $d\ge 2$ the beginning $J_d^0\to J_d^1$ of the coresolution $J_d^\bullet$ equals 
$$\otimes^d\xrightarrow[]{\prod (1-\tau_k)} \bigoplus_{k=0}^{d-2} (\otimes^k)\otimes (\otimes^2)\otimes (\otimes^{d-2-k})\;,$$
where $\tau_k\in\Si_d$ denotes the transposition which exchanges $k+1$ and $k+2$. 
\end{Proposition}
\begin{proof}
The description of the beginning of the coresolution follows from lemma \ref{lm-details-doublebar}. It remains to show that the coresolution is twist compatible. By lemma \ref{lm-puiss-sym-perm-twist-compat}, the maps which define the \emph{algebra} structure of $S^*$ are twist compatible. Now, the differential in the double bar construction $\overline{B}(\overline{B}(A))$ of a CDGA-algebra $A$
is defined using permutations, tensor products and linear combinations of the differential of $A$ and the multiplication of $A$. Thus, by lemmas \ref{lm-puiss-sym-perm-twist-compat} and \ref{lm-twcat-stable}, the differentials of $\overline{B}(\overline{B}(S^*(-)))$ are twist compatible maps.
\end{proof}

\section{Proof of theorem \ref{thm-lifted-classes-cohom-bif}}\label{sec-proof-existence-lifted-classes}
Let $d$ be a positive integer. By proposition \ref{prop-res-Gamma}, the $d$-th
divided powers $\Gamma^d$ admit a twist compatible coresolution $J_d$. By
proposition \ref{prop-twist-compat-res-donne-bicompl}, the totalization of
the associated bicomplex 
$$A(J_d)\,:\, T(J_d^0)_{[1]}(gl)\to T(J_d^1)_{[1]}(gl)
\to \dots \to T(J_d^k)_{[1]}(gl)\to \dots$$
is an $H^*_\mathcal{P}$-acyclic coresolution of $\Gamma^d(gl^{(1)})$, and the
totalization of $H^0_\mathcal{P}(A(J_d))$ computes 
$H^*_\mathcal{P}(\Gamma^d(gl^{(1)})$.
To prove theorem \ref{thm-lifted-classes-cohom-bif}, we exhibit cocycles $z[d]$
representing the lifted classes $c[d]\in H_{\mathcal{P}}^{2d}(\Gamma^d(gl^{(1)}))$
in the bicomplex $H^0_\mathcal{P}(A(J_d))$,
and we prove the relation $\Delta_{(1,\dots,1)\,*}c[d]=c[1]^{\cup d}$ on the 
cochain level.

\subsubsection*{Step 1 : choice of $c[1]$ and $z[1]$}
We first examine the case $d=1$. We may choose a non-zero class $c[1]$ in $H^2_\mathcal{P}(gl^{(1)})\simeq
H^2(GL_{p,\mathbb{k}},\mathfrak{gl}_p^{(1)})$  
(for example we take the `Witt vector class' of \cite[Section 4]{VdK}).

\begin{Notation} We denote by $A_1$ the $p$-coresolution of $gl^{(1)}$ 
obtained by precomposing the $p$-complex $T(S^1)$ by $gl$:
$$A_1:= (T(S^1))(gl) \;.$$ 
\end{Notation}

The bicomplex $A(J_1)$ is concentrated in the first column:
$A(J_1)^{\bullet,
\bullet}=A(J_1)^{0,\bullet}= (A_1)_{[1]}$
and we can choose a cocycle 
$$ z[1]\in \hom_{\mathcal{P}_{p,\mathbb{k}}^p}(\Gamma^p(gl), A_1^p )= \hom_{\mathcal{P}_{p,\mathbb{k}}^p}(\Gamma^p(gl), {A}_{1[1]}^{\;\;2})$$
representing $c[1]$ in $H^0_\mathcal{P}(A(J_1))$ 
(the equality follows from the definition of $-_{[1]}$: the object of degree $2$ of the ordinary complex $A_{1[1]}$ equals the object of degree $p$ of the $p$-complex $A_1$).

\subsubsection*{Step 2 : definition of $z[d]$ and $c[d]$, $d\ge 2$} 
Now we use the cocycle $z[1]$ to build homogeneous cocycles $z[d]$ of bidegree
$(0,2d)$ in the bicomplex $H^0_\mathcal{P}(A(J_d))$, representing the lifted
class $c[d]$ for $d\ge 2$.
\begin{lemme}\label{lm-description-AJ_d}
The first two columns $A(J_d)^{0,\bullet}\to A(J_d)^{1,\bullet}$ of the bicomplex $A(J_d)$ equal:
$$ \underbrace{\left(A_1^{\otimes d}\right)_{[1]}}_{\text{column of index $0$}}\xrightarrow[]{\prod (1-\tau_k)_{[1]}}\underbrace{\bigoplus_{k=0}^{d-2}\left(A_1^{\otimes k}\otimes A_1^{\otimes 2}\otimes A_1^{\otimes d-2-k}\right)_{[1]}}_{\text{column of index $1$}} \;,$$
where the transposition $\tau_k\in\Si_d$ acts by permuting the $(k+1)$-th and the $(k+2)$-th factors of the $p$-complex $A_1^{\otimes d}$.
\end{lemme}
\begin{proof} Use the description of $J_d$ given in proposition \ref{prop-res-Gamma}, and lemma \ref{lm-descrip-sig-mod-Ttens}.
\end{proof}
Let $\Delta_{(p,\dots,p)}:\Gamma^{dp}(gl)\to (\Gamma^p)^{\otimes d}(gl)$ 
be the evaluation of the diagonal $\Gamma^{dp}\to (\Gamma^p)^{\otimes d}$ on 
the bifunctor $gl$. We define:
$$z[d]:=(\underbrace{z[1]\otimes\dots\otimes z[1]}_{\text{$d$ times}})\circ \Delta_{(p,\dots,p)}\in \hom_{\mathcal{P}_{dp,\mathbb{k}}^{dp}}(\Gamma^{dp}(gl), (A_1^p)^{\otimes d}) \;.$$
By lemma \ref{lm-sous-obj}, $(A_1^p)^{\otimes d}$ is a subobject of degree 
$2d$ of the complex $(A_1^{\otimes d})_{[1]}$, so that $z[d]$ is an element of bidegree $(0,2d)$
of the bicomplex $H^0_\mathcal{P}(A(J_d))=\hom(\Gamma^{dp}(gl), A(J_d))$.

\begin{lemme}
$z[d]$ is a cocycle of $\hom_{\mathcal{P}_{dp,\mathbb{k}}^{dp}}
(\Gamma^{dp}(gl), {\mathrm{Tot}} (A(J_d)))$.
\end{lemme}
\begin{proof}
$z[1]$ is a cocycle of degree $2$ in the complex $\hom(\Gamma^p (gl), A_{1[1]}^{\;\;\bullet})$. As a consequence, if $\delta$ denotes the $p$-differential of the $p$-complex $A_1^\bullet$ then $\delta\circ z[1] = 0$. Now the postcomposition by the differential of $(A_1^{\otimes d})_{[1]}$ sends $z[d]$ to the sum
$$\sum_{\ell=1}^d (\underbrace{1\otimes \dots\otimes1\otimes \delta\otimes 1\otimes\dots\otimes 1}_{\text{$\delta$ in $\ell$-th position}})\circ (z[1]^{\otimes d})\circ \Delta_{(p,\dots,p)} \;,$$
and each term of this sum is zero. As a result, the vertical differential of the bicomplex $\hom(\Gamma^{dp}gl, A(J_d)^{\bullet,\bullet})$ sends $z[d]$ to zero. By lemma \ref{lm-description-AJ_d}, the
postcomposition by the horizontal differential of the bicomplex $A(J_d)^{\bullet,\bullet}$ sends $z[d]$ to the sum
$$\sum_{k=0}^{d-2} (1-\tau_k)\circ (z[1]^{\otimes d})\circ \Delta_{(p,\dots,p)}\;.$$
But $\tau_k\circ (z[1]^{\otimes d})\circ \Delta_{(p,\dots,p)}= (z[1]^{\otimes d})\circ \Delta_{(p,\dots,p)}$ so that once again each term of this sum is zero.
Since the horizontal differential and the vertical differential of the
bicomplex $\hom(\Gamma^{dp}(gl), A(J_d)^{\bullet,\bullet})$ both  send $z[d]$
to zero, we deduce that $z[d]$ is a cocycle in the total complex associated to this bicomplex.
\end{proof}
We let $c[d]\in H_{\mathcal{P}}^{2d}(\Gamma^d(gl^{(1)}))$ be the cohomology
class of degree $2d$ represented by the cocycle $z[d]$.

\subsubsection*{Step 3 : proof of the relation $\Delta_{(1,\dots,1)\,*}c[d]=c[1]^{\cup d}$} 
Let $F_1\,A(J_d)$ be the sub-bicomplex of $A(J_d)$ formed by the elements of bidegree $(k,\ell)$ with $k>0$. Thus, $F_1\,A(J_d)$ equals $A(J_d)$ except in the column of index zero where it is null.
The map of complexes
$${\mathrm{Tot}}(A(J_d))\to {\mathrm{Tot}}(A(J_d))/{\mathrm{Tot}}(F_1 A(J_d))=(A_{1}^{\otimes d})_{[1]} $$
is a map of acyclic coresolutions over the diagonal 
$\Delta_{(1,\dots,1)}:\Gamma^{d}(gl^{(1)})\to (gl^{(1)})^{\otimes d}$. 
Thus the cohomology class $\Delta_{(1,\dots,1)\, *}c[d]$ is represented by the
cocycle
$$z[d]=(z[1]^{\otimes d})\circ \Delta_{(p,\dots,p)}\in
\hom_{\mathcal{P}_{dp,\mathbb{k}}^{dp}}(\Gamma^{dp}(gl),(A_1^p)^{\otimes d})
$$ 
in the complex $H^0_\mathcal{P}((A_{1}^{\otimes d})_{[1]})$. On the other
hand, by proposition \ref{prop-cup-acyclic}, the same formula 
$(z[1]^{\otimes d})\circ \Delta_{(p,\dots,p)}$ represents the cup product
$c[1]^{\cup d}$ but this time in another complex, namely
$H^0_\mathcal{P}( (A_{1[1]})^{\otimes d})$. 

So, to finish the proof, we need to compare the two different 
$H^*_{\mathcal{P}}$-acyclic coresolutions $(A_1^{\otimes d})_{[1]}$ 
and $(A_{1[1]})^{\otimes d}$ of $(gl^{(1)})^{\otimes d}$. 
This is achieved by an iterated use of proposition 
\ref{prop-comparaison-prod-tens}: the identity map $\Id_{(gl^{(1)})^{\otimes d}}$ lifts to a map $h^\bullet$ of $H^*_{\mathcal{P}}$-acyclic coresolutions
$$h^\bullet : (A_{1[1]})^{\otimes d} \to (A_1^{\otimes d})_{[1]}\;,$$
such that the restriction of $h^{2d}$ to the subobject $(A_1^p)^{\otimes d}$
is the identity. Thus, postcomposition by $h^\bullet$ sends the cocycle
representing $c[1]^{\cup d}$ to the cocycle representing $\Delta_{(1,\dots,1)\,*}c[d]$. This concludes the proof of theorem \ref{thm-lifted-classes-cohom-bif}.


\begin{thebibliography}{99}
\bibitem{ABW} K.~Akin, D.~Buchsbaum, J.~Weyman, Schur functors and
Schur complexes,  Adv.\ in Math.  44  (1982), 207--278.

\bibitem{CartanSem} H.~Cartan, Constructions multiplicatives. (French) 1958  S\'eminaire Henri Cartan; 7e ann\'ee: 1954/55. Expos\'e no. 4,  8 pp. Secr\'etariat math\'ematique, Paris.

\bibitem{Chalupnik} M.~Cha{\l}upnik, Extensions of strict polynomial functors,
Ann. Sci. \'Ecole Norm. Sup. (4) 38 (2005), no. 5, 773--792.


\bibitem{FF} V.~Franjou, E.~Friedlander, Cohomology of bifunctors, Proc. London Math. Soc. (2008), doi:
10.1112/plms/pdn005.

\bibitem{FFSS} V.~Franjou, E.~Friedlander, A.~Scorichenko, A.~Suslin, General linear and functor cohomology over finite fields,   Ann. of Math. (2)  150  (1999),  no. 2, 663--728.

\bibitem{FLS} V.~Franjou, J.~Lannes, L.~Schwartz, Autour de la cohomologie de Mac Lane des corps finis. (French) [On the Mac Lane cohomology of finite fields]  Invent. Math.  115  (1994),  no. 3, 513--538.

\bibitem{FS} E.~Friedlander, A.~Suslin, Cohomology of finite group schemes over a field,
 Invent. Math. 127 (1997), 209--270.
 
\bibitem{Tohoku} A.~Grothendieck, Sur quelques points d'alg\`ebre homologique. (French)  T�oku Math. J. (2)  9  (1957), 119--221.

\bibitem{Haboush} W.~Haboush, Reductive groups are geometrically reductive,  Ann. of Math. (2)  102  (1975), no. 1, 67--83.
 

\bibitem{Kap} M.M.~Kapranov, On the q-analog of homological algebra, \texttt{arXiv:q-alg/9611005} (1996).

\bibitem{KW} C.~Kassel, M.~Wambst, Alg\`ebre homologique des $N$-complexes et homologie de Hochschild aux racines de l'unit\'e. (French) [Homology algebra of $N$-complexes and Hochschild homology at the roots of unity]  Publ. Res. Inst. Math. Sci.  34  (1998),  no. 2, 91--114.

\bibitem{ML} S.~Mac~Lane, Homology. Reprint of the 1975 edition. Classics in Mathematics. Springer-Verlag, Berlin, 1995. x+422 pp. ISBN: 3-540-58662-8.

\bibitem{Nagata} M.~Nagata,  Invariants of a group in an affine ring, J.\ Math.\ Kyoto Univ. 3 (1963/1964), 369--377.

\bibitem{Real} P.~Real, Homological perturbation theory and associativity.  Homology Homotopy Appl.  2  (2000), 51--88 (electronic).

\bibitem{SvdK}V.~Srinivas and W.~van der Kallen, Finite Schur filtration dimension for modules over an algebra with Schur
filtration, to appear in Transform. Groups


\bibitem{Totaro} B. Totaro, Projective resolutions of representations of ${\rm GL}(n)$.  J. Reine Angew. Math.  482  (1997), 1--13.

\bibitem{Touze} A. Touz\'e, Ph.D. thesis (2008)

\bibitem{TvdK} A. Touz\'e, W.~van der Kallen, Bifunctor cohomology and
cohomological finite generation for reductive groups, 
{\tt arXiv:0809.1014}. 

\bibitem{Troesch} A.~Troesch, Une r\'esolution injective des puissances sym\'etriques tordues. (French) [Injective resolution of twisted symmetric powers]  Ann. Inst. Fourier (Grenoble)  55  (2005),  no. 5, 1587--1634. 

\bibitem{VdK} W.~van der Kallen, Cohomology with Grosshans graded coefficients, In: Invariant Theory in All Characteristics, Edited by: H. E. A. Eddy Campbell and David L. Wehlau, CRM Proceedings and Lecture Notes, Volume 35 (2004) 127-138, Amer. Math. Soc., Providence, RI, 2004.

\bibitem{VdK2} W.~van der Kallen, A reductive group with finitely generated cohomology algebras.
In: Algebraic Groups and Homogeneous Spaces, Mumbai 2004, Edited by: Vikram B. Mehta. Narosa, 2007. ISBN:  978-81-7319-802-1
\end{thebibliography}
\end{document}